\documentclass[12pt]{amsart}

\usepackage[margin=1.5in]{geometry} 

\usepackage[T1]{fontenc}
\usepackage[utf8]{inputenc}
\usepackage[pdftex]{graphicx}
\usepackage{tikz-cd}
\usepackage{blindtext}

\usepackage{hyperref,xcolor}
\usepackage{amsmath,amsfonts,amssymb,amsthm,color,mathtools,enumitem}
\usepackage{epstopdf}
\usepackage{polynom}
\usepackage{babel}
\usepackage{mathrsfs}
\usepackage{chngcntr}
\usepackage{verbatim}
\usepackage{hyphenat}

\hypersetup{
	colorlinks=true,
	linkcolor=blue,
	citecolor=cyan,
	urlcolor=cyan}

\newtheorem{theorem}{Theorem}[section]
\newtheorem{defn}[theorem]{Definition}
\newtheorem{prop}[theorem]{Proposition}
\newtheorem{cor}[theorem]{Corollary}
\newtheorem{lemma}[theorem]{Lemma}
\newtheorem{example}[theorem]{Example}
\newtheorem{remark}[theorem]{Remark}

\newtheorem{result}{Theorem}

\newtheorem{resultcor}[result]{Corollary}

\numberwithin{equation}{section}

\newcommand{\into}{\hookrightarrow}

\newcommand{\act}{\curvearrowright}

\newcommand{\ee}{\varepsilon}
\newcommand{\sm}{\setminus}

\newcommand{\dist}{\text{dist}}

\newcommand{\id}{\text{id}}

\newcommand{\Ad}{\text{Ad}}

\newcommand{\II}{\text{II}}

\DeclareMathOperator{\End}{End}
\DeclareMathOperator{\Aut}{Aut}

\newcommand{\dlim}{\underset{\to}{\lim}}

\newcommand{\ov}{\overline}
\newcommand{\bN}{\mathbb{N}}
\newcommand{\bR}{\mathbb{R}}
\newcommand{\bZ}{\mathbb{Z}}
\newcommand{\bQ}{\mathbb{Q}}
\newcommand{\bF}{\mathbb{F}}

\newcommand{\bC}{\mathbb{C}}

\newcommand{\cM}{\mathcal{M}}
\newcommand{\cN}{\mathcal{N}}

\newcommand{\cO}{\mathcal{O}}

\newcommand{\cL}{\mathcal{L}}
\newcommand{\cF}{\mathcal{F}}

\newcommand{\cD}{\mathcal{D}}

\newcommand{\cZ}{\mathcal{Z}}

\newcommand{\cG}{\mathcal{G}}
\newcommand{\cR}{\mathcal{R}}
\newcommand{\cQ}{\mathcal{Q}}

\begin{document}
\title{Tensorially absorbing inclusions of C*-algebras}
\author{Pawel Sarkowicz}
\email{\href{mailto:psark007@uottawa.ca}{psark007@uottawa.ca}}
\address{Department of Mathematics and Statistics, University of Ottawa, 75 Laurier Ave. East, Ottawa, ON, K1N 6N5 Canada}

  \begin{abstract}
    When $\mathcal{D}$ is strongly self-absorbing we say an inclusion $B \subseteq A$ is $\mathcal{D}$-stable if it is isomorphic to the inclusion $B \otimes \mathcal{D} \subseteq A \otimes \mathcal{D}$. We give ultrapower characterizations and show that if a unital inclusion is $\mathcal{D}$-stable, then $\mathcal{D}$-stability can be exhibited for countably many intermediate C*-algebras concurrently. We show that such embeddings between $\mathcal{D}$-stable C*-algebras are point-norm dense in the set of all embeddings, and that every embedding between $\mathcal{D}$-stable C*-algebras is approximately unitarily equivalent to a $\mathcal{D}$-stable embedding. Examples are provided. 
  \end{abstract}

  \maketitle
  \tableofcontents

\section{Introduction}
  The study of inclusions of C*-algebras has been of recent interest. There is no short supply of research concerning inclusions relating to non-commutative dynamics \cite{Popa00,Izumi02,CameronSmith19,OsakaTeruya18,EchterhoffRordam21}, as well as inclusions of simple C*-algebras \cite{Rordam21}. We discuss inclusions from the lens of tensorially absorbing a strongly self-absorbing C*-algebra $\cD$ \cite{TomsWinter07}.

  When speaking of tensorial absorption with a strongly self-absorbing C*-algebra, central sequences play an imperative role akin to McDuff's character\hyp{}ization of when a $\II_1$ von Neumann algebra absorbs the unique hyperfinite $\II_1$ factor $\cR$ \cite{McDuff69}. Central sequences have been studied since the inception of operator algebras as Murray and von Neumann used them to exhibit non-isomorphic $\II_1$ factors by showing that $\cL(\bF_2)$ does not have property $\Gamma$ \cite{MvNIV}. They were also used in Connes' theorem concerning the uniqueness of $\cR$ \cite{Connes76}, and the classification of automorphisms on hyperfinite factors \cite{Connes75,Connes76}. In \cite{Bisch90,Bisch94}, Bisch considered the central sequence algebra $\cN^\omega \cap \cM'$ associated to an (irreducible) inclusion of $\II_1$ factors $\cN \subseteq \cM$ and characterized when there was an isomorphism $\Phi: \cM \simeq \cM \ov{\otimes} \cR$ such that $\Phi(\cN) = \cN \ov{\otimes} \cR$ in terms of there being non-commuting sequences in $\cN$ which asymptotically commute with the larger von Neumann algebra $\cM$ (in the $\|\cdot\|_2$-norm). As pointed out by Izumi \cite{Izumi04}, there are similar central characterizations for unital inclusions of separable C*-algebras which tensorially absorb a strongly self-absorbing C*-algebra $\cD$ (it was at least pointed out for $\cD$ being one of $M_{n^{\infty}},\cO_2,\cO_{\infty}$).

  For a strongly self-absorbing C*-algebra $\cD$ \cite[Definition 1.3(iv)]{TomsWinter07}, we study $\cD$-stable
  inclusions (see Section \ref{section:dse-embeddings} for detailed definitions), analogous to Bisch's notion for an (irreducible) inclusion of $\II_1$ factors \cite{Bisch90}. We say that an inclusion $B \subseteq A$ is $\cD$-stable if there is an isomorphism $A \simeq A \otimes \cD$ such that
\begin{equation}
  \begin{tikzcd}
A \arrow[r, "\simeq"]                 & A \otimes \cD                 \\
B \arrow[r, "\simeq"] \arrow[u, hook, "\iota"] & B \otimes \cD \arrow[u, hook, "\iota \otimes \id_\cD"']
\end{tikzcd}
\end{equation}
commutes.

We study such inclusions systematically, discussing central sequence char\hyp{}acterizations, permanence properties, and giving examples towards the end. We list some key findings here.  The first is that $\cD$-stable inclusions exist between $\cD$-stable C*-algebras if there is any inclusion, and that the set of $\cD$-stable inclusions is quite large. Moreover, as far as classification of embeddings up to approximate unitary equivalence (in particular by $K$-theory and traces), $\cD$-stable embeddings are all that matter. 

\begin{result}[Proposition \ref{allembeddingsareue}, Corollary \ref{pointnormdensity}]
    Let $A,B$ be unital, separable, $\cD$-stable C*-algebras.
      \begin{enumerate}
        \item The set of $\cD$-stable embeddings $B \into A$ is point-norm dense in the set of all embeddings $B \into A$.
        \item Every embedding $B \into A$ is approximately unitarily equivalent to a $\cD$-stable embedding. 
      \end{enumerate}
  \end{result}

  We note that this set is however not everything. We provide examples of non-$\cD$-stable inclusions of $\cD$-stable C*-algebras, namely by fitting a C*-algebra with perforated Cuntz semigroup or with higher stable rank (in partic\hyp{}ular non-$\cZ$-stable C*-algebras) in between two $\cD$-stable C*-algebras. The second useful tool is that a $\cD$-stable inclusion allows one to find an appropriate isomorphism witnessing $\cD$-stability of countably many intermediate subalge\hyp{}bras at once.

\begin{result}[Theorem \ref{countablymanysubalgebras}]
    Let $B \subseteq A$ be a unital, $\cD$-stable inclusion of separable C*-algebras. If $(C_n)$ is a sequence of C*-algebras such that $B \subseteq C_n \subseteq A$ unitally for all $n$, then there exists an isomorphism $\Phi: A \simeq A \otimes \cD$ such that
      \begin{enumerate}
        \item $\Phi(B) = B \otimes \cD$ and 
        \item $\Phi(C_n) = C_n \otimes \cD$ for all $n \in \bN$. 
      \end{enumerate}
  \end{result}
This is not a trivial condition, as it is not true that any such isomorphism sends every intermediate C*-algebra to its tensor product with $\cD$ (see Example \ref{noteveryintermediate}). In fact, one can always find an intermediate C*-algebra $C$ between $B$ and  $A$ and an isomorphism $A \simeq A \otimes \cD$ sending $B$ to $B \otimes \cD$ which does not send $C$ to $C \otimes \cD$ (although, of course, we will still have $C \simeq C \otimes \cD$). 

The above result, together with the Galois correspondence of Izumi \cite{Izumi02}, allows us to get a result similar to the main theorem of \cite{AGJP22}. There they prove that if $G \act^\alpha A$ is an action of a finite group with the weak tracial Rokhlin property on a C*-algebra $A$ with sufficient regularity conditions, then every C*-algebra between $A^\alpha \subseteq A$ and $A \subseteq A \rtimes_\alpha G$ is $\cZ$-stable. Assuming we have a unital C*-algebra with the same regularity conditions, we show that we can witness $\cZ$-stability of all such intermediate C*-algebras concurrently.

  \begin{resultcor}[Corollary \ref{weaktracialRokhlinintermediate}]
            Let $A$ be a unital, simple, separable, nuclear $\cZ$-stable C*-algebra and $G \act^\alpha A$ be an action of a finite group with the weak tracial Rokhlin property. There exists an isomorphism $\Phi: A \rtimes_\alpha G \simeq (A \rtimes_\alpha G) \otimes \cZ$ such that whenever $C$ is a unital C*-algebra satisfying either
              \begin{enumerate}
                \item $A^\alpha \subseteq C \subseteq A$ or
                \item $A \subseteq C \subseteq A \rtimes_\alpha G$,
              \end{enumerate}
              we have $\Phi(C) = C \otimes \cZ$. 
          \end{resultcor}

  This paper is structured as follows. We discuss various local properties in Section \ref{section:dse-localembeddings}, and then formalize the notion of a $\cD$-stable embedding in Section \ref{section:dse-embeddings}, examining several properties and consequences. In Section \ref{section:dse-crossed-products} we show how several examples arising from non-commutative dynamical systems fit into the framework of $\cD$-stable inclusions. We finish with several examples in Section \ref{section:dse-examples}.



  \addtocontents{toc}{\protect\setcounter{tocdepth}{0}}
  \section*{Acknowledgements}
  Many thanks to my supervisors Thierry Giordano and Aaron Tikuisis for many helpful discussions, as well as to Eusebio Gardella for many helpful comments.  

  \addtocontents{toc}{\protect\setcounter{tocdepth}{1}}

    \section{Preliminaries}\label{preliminaries}

  \subsection{Notation} We will use capital letters $A,B,C,D$ to denote C*-algebras and usually a calligraphic $\cD$ to denote a strongly self-absorbing C*-algebra. Generally small letters $a,b,c,d,\dots,x,y,z$ will denote operators in C*-algebras. $A_+$ will denote cone of positive elements in a C*-algebra $A$. If $\ee > 0$ and $a,b$ are elements in a C*-algebra, we will write
\begin{equation}
   a \approx_\ee b
\end{equation}
  to mean that $\|a - b\| < \ee$. This will make some approximations more legible. 

  The symbol $\otimes$ will denote the minimal tensor product of C*-algebras, while $\odot$ will mean the algebraic tensor product. We use the minimal tensor product throughout, and it is common for us to deal with nuclear C*-algebras so there should not be any ambiguity. The symbol $\ov{\otimes}$ will denote the von Neumann tensor product. 

  We will denote by $M_n$ the C*-algebra of $n\times n$ matrices, and $M_{n^{\infty}}$ the uniformly hyperfinite (UHF) C*-algebra associated to the supernatural number $n^{\infty}$. We will write $\cQ$ for the universal UHF algebra $\cQ = \bigotimes_{n \in \bN} M_n$. 

  By $G \act^\alpha A$, we will mean that the (discrete) group $G$ acts on $A$ by automorphisms, i.e., $\alpha: G \to \Aut(A)$ is a homomorphism. $A \rtimes_{r,\alpha}G$ will denote the reduced crossed product, which we will just write as $A \rtimes_\alpha G$ if it is clear from context that the group is amenable and $A$ is nuclear (e.g., if $G$ is finite). We will denote by $A^\alpha$ the fixed point subalgebra of the action (or $A^G$ if the action is clear from context).

\subsection{Ultrapowers, central sequences and central sequence algebras}\label{section:ultrapowers}

Fix a free ultrafilter $\omega \in \beta \bN$. Throughout we will use ultrapowers to describe asymptotic behaviour, rather than sequences algebras. This comes down to a matter of taste and one can swap between the two if one so desires, as we will provide local characterizations as well. This also means that all of what we do will be independent of the specific ultrafilter $\omega$.

For a C*-algebra $A$, the ultrapower of $A$ is the C*-algebra
  \begin{equation}
 A_\omega := \ell^{\infty}(A)/c_{0,\omega}(A), 
  \end{equation}
where $c_{0,\omega} := \{(a_n) \in \ell^{\infty}(A) \mid \lim_{n \to \omega} \|a_n\| = 0\}$ is the ideal of $\omega$-null sequences. We can embed $A$ into $A_\omega$ canonically by means of constant sequences: we identify $a \in A$ with the equivalence class of the constant sequence $(a)_n$.

To ease notation, we will usually write elements of $A_\omega$ as sequences $(a_n)$, keeping in mind that these are equivalence classes without explicitly stating it every time. We note that the norm on $A_\omega$ is given by $\|(a_n)\| = \lim_{n \to \omega} \|a_n\|$. 

Kirchberg's $\ee$-test (\cite{Kirchberg06}, Lemma A.1) is essentially the operator algebraists' {\L}o{\' s}' theorem without having to turn to (continuous) model theory. Heuristically, it says that if certain things can be done approximately in an ultrapower, then certain things can be done exactly in an ultrapower.

\begin{lemma}[Kirchberg's $\ee$-test]\label{lemma:Kirchberg-eptest}
Let $(X_n)_n$ be a sequence of sets and suppose that for each $n$, there is a sequence $(f_n^{(k)})$ of functions $f_n^{(k)}: X_n \to [0,\infty)$. For $k \in \bN$, let
\begin{equation}
f_\omega^k(s_1,s_2,\dots) := \lim_{n \to \omega} f_n^{(k)}(s_n). 
\end{equation}
Suppose that for every $m \in \bN$ and $\ee > 0$, there is $s \in \prod_n X_n$ with $f_\omega^{(k)}(s) < \ee$ for $k=1,\dots,m$. Then there exists $t \in \prod_nX_n$ with $f_\omega^{(k)}(t) = 0$ for all $k \in \bN$. 
\end{lemma}

The above is useful, although if one so wishes, one can usually construct exact objects from approximate objects by using standard diagonalization arguments (under some separability assumptions).
These sorts of arguments work in both the ultrapower setting and the sequence algebra setting.

Finally, if $\alpha \in \Aut(A)$ is an automorphism, there is an induced automorphism on $A_\omega$, which we will denote by $\alpha_\omega$ given by
\begin{equation}
\alpha_\omega((a_n)) := (\alpha(a_n)). 
\end{equation}

    \subsection{Central sequences and central sequence subalgebras} For a unital C*-algebra $A$, the C*-algebra of $\omega$-central sequences is
\begin{equation}
    A_\omega \cap A' = \{x \in A_\omega \mid [x,a] = 0 \text{ for all } a \in A\},
\end{equation}
    where we are identifying $A \subseteq A_\omega$ with the constant sequences. If $B \subseteq A$ is a unital C*-subalgebra and $S \subseteq A_\omega$ is a subset, we can associate the relative commutant of $S$ in $B_\omega$:
\begin{equation}
    B_\omega \cap S' = \{b \in B_\omega \mid [b,s] = 0 \text{ for all } s \in S\}. 
\end{equation}
    Of particular interest will be when $S = A$, and $B \subseteq A$ is a unital inclusion of separable C*-algebras.

\subsection{Strongly self-absorbing C*-algebras}\label{section:ssa-C*-algebras}

A unital separable C*-algebra $\cD$ is strongly self-absorbing if $\cD \not\simeq \bC$ and there is an isomorphism $\phi: \cD \to \cD \otimes \cD$ which is approximately unitarily equivalent to the first factor embedding $d \mapsto d \otimes 1_\cD$ (see \cite{TomsWinter07}). All known strongly self-absorbing C*-algebras are: the Jiang-Su algebra $\cZ$ \cite{JiangSu99}, the Cuntz algebras $\cO_2$ and $\cO_{\infty}$ \cite{Cuntz77}, UHF algebras of infinite type, and $\cO_{\infty}$ tensor a UHF algebra of infinite type. Strongly self-absorbing C*-algebras have approximately inner flip, and therefore there are $K$-theoretic restrictions on strongly self-absorbing C*-alge\hyp{}bras -- see \cite{Tikuisis16,EndersSchemaitatTikuisis23}.

Tensorial absorption with strongly self-absorbing C*-algebras gives rise to many regular properties, for examples in terms of $K$-theory, traces, and the Cuntz semigroup \cite{JiangSu99,Rordam91,Rordam92,Rordam04}. Of paramount interest is the  Jiang-Su algebra $\cZ$. A cumulation of work has successfully classified all (unital) separable, simple, nuclear, infinite-dimensional, $\cZ$-stable C*-algebras satisfying the Universal Coefficient Theorem (UCT) of Rosenberg and Schochet \cite{RosenbergSchochet87} by means of $K$-theory and traces. We describe how one might work with $\cZ$-stability in terms of its standard building blocks. Recall that, for $n,m \geq 2$, the dimension drop algebras are
  \begin{equation}
\cZ_{n,m} := \{f \in C([0,1],M_n \otimes M_m) \mid f(0) \in M_n \otimes 1_{M_m}, f(1) \in 1_{M_n} \otimes M_m\}. 
  \end{equation}
Such an algebra is a called a prime dimension drop algebra when $n$ and $m$ are coprime. The Jiang-Su algebra $\cZ$ is the unique separable simple C*-algebra with unique tracial state which is an inductive limit of prime dimension drop algebras with unital connecting maps \cite{JiangSu99} (in fact, the dimension drop algebras can be chosen to have the form $\cZ_{n,n+1}$).It is $KK$-equivalent to $\bC$ and $\cZ$-stability is a necessary condition for $K$-theoretic classification.

By \cite[Proposition 5.1]{RordamWinter10} (or \cite[Proposition 2.1]{Sato10} for our desired formulation), $\cZ_{n,n+1}$ is the universal C*-algebra generated by elements $c_1,\dots,c_n$ and $s$ such that
\begin{itemize}
\item $c_1 \geq 0$;
\item $c_ic_j^* = \delta_{ij}c_1^2$;
\item $s^*s + \sum_{i=1}^n c_i^*c_i = 1$;
\item $c_1s = s$.
\end{itemize}
If there are uniformly tracially large (in the sense of \cite[Definition 2.2]{TomsWhiteWinter15}) order zero\footnote{order zero meaning orthogonality preserving: $\phi: A \to B$ is c.p.c. order zero if it is c.p.c. and $\phi(a)\phi(b) = 0$ whenever $ab = 0$.} c.p.c. maps $M_n \to A_\omega \cap A'$, these give rise to elements $c_1,\dots,c_n \in A_\omega \cap A'$ with $c_1 \geq 0$ and $c_ic_j^* = \delta_{ij}c_1^2$, along with certain tracial information. If $A$ has strict comparison, Matui and Sato used this tracial information to show that $A$ has property (SI) \cite{MatuiSato12}, from which one can get an element $s \in A_\omega \cap A'$ such that $s^*s + \sum_{i=1}^n c_i^*c_i = 1$ and $c_1s = s$. This gives a *-homomorphism $\cZ_{n,n+1} \to A_\omega \cap A'$, which if can be done for each $n \in \bN$, is enough to conclude that $\cZ \into A_\omega \cap A'$ unitally and hence $A \simeq A \otimes \cZ$. In fact, it suffices to show that $\cZ_{2,3} \into A_\omega \cap A'$ (or $\cZ_{n,n+1}$ for some $n \geq 2$), see \cite[Theorem 3.4(ii)]{RordamWinter10} and \cite[Theorem 5.15]{Schemaitat22}. 

\section{Approximately central approximate embeddings}\label{section:dse-localembeddings}

  Here we formalize some results on approximate embeddings. When $B \subseteq A$ is a unital inclusion of separable C*-algebras, this will yield local characterizations of nuclear subalgebras of $B_\omega \cap A'$, as defined in (\ref{eq:central-sequence-subalgebra}).

    \begin{defn}
      Let $B \subseteq A$ be a unital inclusion of C*-algebras and let $D$ be a unital, simple, nuclear C*-algebra. Let $\cF \subseteq D, \cG \subseteq A$ be finite sets and $\ee > 0$. We say that a u.c.p. map $\phi: D \to B$ is an $(\cF,\ee)$-approximate embedding if
        \begin{enumerate}
          \item $\phi(cd) \approx_\ee \phi(c)\phi(d)$ for all $c,d \in \cF$.
        \end{enumerate}
        If $\phi$ additionally satisfies
          \begin{enumerate}
            \item[2.] $[\phi(c),a] \approx_\ee 0$ for all $c \in \cF$ and $a \in \cG$,
          \end{enumerate}
          then we say that $\phi$ is an $(\cF,\ee,\cG)$-approximately central approximate embedding. 
    \end{defn}

    We will usually write that $\phi$ is a $(\cF,\ee)$-embedding or $(\cF,\ee,\cG)$-embedding to mean that $\phi$ is an $(\cF,\ee)$-approximate embedding or $(\cF,\ee,\cG)$-approximately central approximate embedding respectively. 

  \begin{remark}
    One can make a similar definition to the above if $D$ is not simple or nuclear (or even unital). The aim is to discuss subalgebras of $B_\omega \cap A'$, and if $D \into B_\omega \cap A'$ is nuclear, then one can use Choi-Effros to lift the embedding to a sequence of u.c.p. maps which are approximatey isometric, approximately multiplicative, and approximately commute with finite subsets of $A$. If $D$ is simple, the approximate isometry condition follows for free since the embedding $D \into B_\omega \cap A'$ must be isometric.
    
If we loosen the simple and nuclear assumptions on $D$, we can still speak of bounded linear maps $\phi:D \to B$ (no longer necessarily u.c.p.) which are approximately isometric, approximately multiplicative, approximately adjoint-preserving, and approximately commute with a finite prescribed subset of $A$. This will allow one to discuss general subalgebras of $B_\omega \cap A'$. As we will only be interested in strongly self-absorbing subalgebras of $B_\omega \cap A'$, which are unital, separable, simple, and nuclear (see \cite[Section 1.6]{TomsWinter07}), we restrict ourselves to u.c.p. maps from a unital, simple, nuclear C*-algebras which are approximately multiplicative and approximately commute with finite subsets of $A$.

Most of the work in this section can be done without assumptions of simplicity and nuclearity. 
  \end{remark}

    \begin{lemma}\label{lem:local-equiv-global}
      Suppose that $A,B,D$ are unital C*-algebras with $B$ separable and $D$ simple, separable and nuclear. Let $S \subseteq A$ be a separable subset. There are $(\cF,\ee,\cG)$-approximately central approximate embeddings $D \to B$ for all $\cF \subseteq D, \cG \subseteq S$ and $\ee > 0$ if and only if there is a unital embedding $D \into B_\omega \cap S'$.
        \begin{proof}
          Let $(F_n)$ be an increasing sequence of finite subsets of $D$ with dense union and let $(G_n)$ be an increasing sequence of finite subsets of $S$ with dense union. Let $\phi_n: D \to B$ be  $(F_n,\frac{1}{n},G_n)$-approximately central approximate embeddings. Let $\pi: \ell^{\infty}(B) \to B$ denote the quotient map and set
          \begin{equation}
                     \psi := \pi \circ (\phi_n): D \to B_\omega
\end{equation} 
              which is a unital embedding such that $[\psi(d),a] = 0$ for  all $d \in D$ and $a \in S$. 

          For the other direction, suppose that $\psi: D \to B_\omega \cap S'$ is a unital embedding, $\cF \subseteq D,\cG \subseteq S$ are finite and $\ee >0$. By the Choi-Effros lifting theorem (see, for example, \cite[Theorem C.3]{BrownOzawa}) there is a u.c.p. lift $\tilde{\psi} = (\tilde{\psi}_n): D \to \ell^{\infty}(B)$ such that
            \begin{itemize}
              \item $\|\tilde{\psi}_n(cd) - \tilde{\psi}_n(c)\tilde{\psi}_n(d)\| \to^{n \to \omega} 0$,
              \item $\|[\tilde{\psi}_n(d),a]\| \to^{n \to \omega} 0$ 
            \end{itemize}
            for all $c,d \in D$ and $a \in A$. Take $n$ large enough and set $\phi = \psi_n$, so that $\phi$ will be a $(\cF,\ee,\cG)$-approximately central approximate embedding. 
        \end{proof}
    \end{lemma}

      \begin{cor}
        Let $A,B,D$ be unital C*-algebras with $B,D$ separable, simple and nuclear. Suppose that there are unital embeddings $\phi: D \to B_\omega$ and $\psi: B \to A_\omega$.  Then there is a unital embedding $\xi: D \into A_\omega$. If $S \subseteq A_\omega$ is a separable subset with $\psi(B) \subseteq A_\omega \cap S'$, then $\xi$ can be chosen with $\xi(C) \subseteq A_\omega \cap S'$. 
          \begin{proof}
            Let $\cF \subseteq D$ be finite and $\ee > 0$. Let $L := \max\{\max_{d \in \cF}\|d\|,1\}$. By the above lemma, there is an $(\cF,\frac{\ee}{2L})$-approximate embedding $\phi: D \to B$, so let $\cF' = \phi(\cF)$. Now there is an $(\cF',\frac{\ee}{2L})$-approximate embedding $\psi: B \to A$. An easy calculation shows that $\psi \circ \phi: D \to A$ is an approximate $(\cF,\ee)$-embedding. 

            Appending the condition that $\psi: B \to A_\omega  \cap S'$, then, for any finite subset $\cG \subseteq S$, we can take $\psi: B \to A$ to be a $(\cF',\frac{\ee}{2L},\cG)$-approximately central approximate embedding. This gives that $\psi \circ \phi: D \to A$ is be a $(\cF,\ee,\cG)$-approximately central approximate embedding. 
          \end{proof}
      \end{cor}

        \begin{cor}
          Let  $D$ be a C*-algebra and $B \subseteq A$ be a unital inclusion of separable C*-algebras such that $B$ and $D$ are unital, separable, simple and nuclear. Suppose that there is an embedding $\pi: A \into A_\omega \cap A'$ with $\pi(B) \subseteq B_\omega \cap A'$. If $D \into B_\omega$ unitally, then $D \into B_\omega \cap A'$ unitally.
            \begin{proof}
              As $D \into B_\omega$ and $B \into B_\omega \cap A' \subseteq A_\omega \cap A'$, the above yields $D \into B_\omega \cap A'$. 
            \end{proof}
        \end{cor}

The following is useful for discussing $\cD$-stability for some inclusions of fixed point subalgebras by certain automorphisms on UHF algebras.  In particular, the following will work for automorphisms on UHF algebras of product-type, as well as tensor permutations (of finite tensor powers of UHF algebras). 

\begin{cor}\label{autfixedpoint}
  Let $A = \bigotimes_\bN B$ be an infinite tensor product of a unital, separable, nuclear C*-algebra $B$ and let $D$ be unital, separable, simple, and nuclear. Let $\lambda \in \End(A)$ be the Bernoulli shift $\lambda(a) = 1 \otimes a$. If $\sigma \in \Aut(A)$ is such that $\lambda \circ \sigma  = \sigma \circ \lambda$, and $D \into (A^\sigma)_\omega$ unitally, then $D \into (A^\sigma)_\omega \cap A'$ unitally. 
              \begin{proof}
                Note that $\pi = (\lambda^n)$ induces an embedding $A \into A_\omega \cap A'$. We just need to show that $\pi(A^\sigma) \subseteq (A^\sigma)_\omega \cap A'$. The hypothesis gives that $\lambda^n \circ \sigma = \sigma \circ \lambda^n$ for all $n$, hence $\pi(A^\sigma) \subseteq (A^\sigma)_\omega \cap A'$. The result now follows from the above. 
              \end{proof}
          \end{cor}

  We note that if we have approximately central approximate embeddings $D \to B \subseteq A$, then we can also find approximately central approximate embedding $D \to u^*Bu \subseteq A$ for any $u \in U(A)$. In the separable setting, this just means $D \into B_\omega \cap A'$ implies that $D \into u^*B_\omega u \cap A'$ for any $u \in U(A)$. 

  \begin{lemma}\label{aduaapproximateembeddings}
    Let $B \subseteq A$ be a unital inclusion of C*-algebras and let $D$ be a unital, separable, simple, nuclear C*-algebra. Let $u \in U(A)$. If there are $(\cF,\ee,\cG)$-approximately central approximate embeddings $D \to B$ for all $\cF \subseteq D, \cG \subseteq A$ finite subsets and $\ee > 0$, then there are $(\cF,\ee,\cG)$-approximately central approximate embeddings $D \to u^*Bu \subseteq A$ for all $\cF,\ee,\cG$.
      \begin{proof}
        Let $\cF \subseteq D, \cG \subseteq A$ be finite and $\ee > 0$. Let $L = \max\{1,\max_{d \in \cF}\|d\|\}$ and $\phi: D \to B$ be a $(\cF,\frac{\ee}{3L},\cG \cup \{u\})$-approximately central approximate embedding. Then $\psi = \Ad_u \circ \phi: D \to u^*Bu$ will be an $(\cF,\ee,\cG)$-embedding. 
      \end{proof}
  \end{lemma}

  We can also discuss existence of approximately central approximate embed\hyp{}dings in inductive limits (with injective connecting maps). This is an adaptation of \cite[Proposition 2.2]{TomsWinter08} to our setting. 

  \begin{prop}\label{approximateembeddinginductivelimit}
      Suppose that we have increasing sequences $(B_n)$ and $(A_n)$ of C*-algebras such that $B_n \subseteq A_n$ are unital inclusions.
If $B = \ov{\cup_n B_n}, A = \ov{\cup_n A_n}$, and $D = \ov{\cup_nD_n}$ where $(D_n)$ is an increasing sequence of unital, separable, simple, nuclear C*-algebras and there are $(\cF,\ee,\cG)$-embeddings $D_n \to B_n \subseteq A_n$ whenever $n \in \bN, \cF \subseteq D_n, \cG \subseteq A_n$ are finite and $\ee > 0$, then there are $(\cF,\ee,\cG)$-embeddings $D \to B \subseteq A$ for all $\cF \subseteq D,\cG \subseteq A$ finite and $\ee > 0$. 
  \begin{proof}
    Let $\cF \subseteq \cD$ and $\cG \subseteq A$ be finite sets and $\ee > 0$. Let
\begin{equation}
    L := \max\{1,\max_{d \in \cF}\|d\|,\max_{a \in \cG}\|a\|\}
\end{equation}
     and set $\delta := \frac{\ee}{6L+5}$. Without loss of generality assume that $\ee < 1$. Label $\cF = \{d_1,\dots,d_p\}$ and $\cG = \{a_1,\dots,a_q\}$ and find $N$ large enough so that there are $d_1',\dots,d_p' \in D_N$ and $a_1',\dots,a_q' \in A_N$ with $d_i' \approx_\delta d_i$ and $a_j' \approx_\delta a_j$. Let $\cF':=\{d_1',\dots,d_p'\}, \cG':=\{a_1',\dots,a_q'\}$ and let $\phi: D_N \to B_N \subseteq A_N$ be an $(\cF',\delta,\cG')$-embedding. As $D_N$ is nuclear, there are $k \in \bN$ and u.c.p. maps $\rho: D_N \to M_k$ and $\eta: M_k \to B_N$ such that $\eta \circ \rho(d_i') \approx_\delta \phi(d_i')$ and $\eta \circ \rho(d_i'd_j') \approx_\delta \phi(d_i'd_j')$. Use Arveson's extension theorem (see \cite[Section 1.6]{BrownOzawa}) to extend $\rho$ to a u.c.p. map $\tilde{\rho}: D \to M_k$ and let $\psi := \eta \circ \tilde{\rho}: D \to B_N$. As $B_N \subseteq B$, we can think of $\psi$ as a map $\psi: D \to B$. Now for $i=1,\dots,p$, we have
\begin{equation}
\begin{split}
        \psi(d_id_j) &\approx_{(2L+1)\delta}\psi(d_i'd_j') \\
                     &= \eta \circ \rho(d_i'd_j') \\
                     &\approx_\delta \phi(d_i'd_j') \\
                     &\approx_\delta \phi(d_i')\phi(d_j') \\
                     &\approx_{2L\delta} \eta \circ \rho(d_i')\eta \circ \rho(d_j') \\
                     &= \psi(d_i')\psi(d_j') \\
                     &\approx_{(2L+1)\delta}\psi(d_i)\psi(d_j).
\end{split}
\end{equation}
      Thus $\psi(d_id_j) \approx_{(4+6L)\delta} \psi(d_i)\psi(d_j)$, and as $(4 + 6L)\delta \leq (6L + 5)\delta = \ee$, this implies that $\psi(d_id_j) \approx_\ee \psi(d_i)\psi(d_j)$. For approximate commutation with $\cG$, we make use of the following two approximations: for $a,a',a'',b,b'$ elements in a C*-algebra,
      \begin{equation}
      \begin{split}
      &\|[a,b]\|  \leq (\|a\| + \|a'\|)\|b - b'\| + (\|b\| + \|b'\|)\|a - a'\| + \|[a',b']\|, \\
      &\|[a',b']\| \leq 2\|b'\|\|a' - a''\| + \|[a'',b']\|.
      \end{split}
      \end{equation}
        Note that for $a = \psi(d_i),a' = \psi(d_i'),a'' = \phi(d_i'), b = a_j, b' = a_j'$, we have that $\|a\|,\|b\| \leq L + 1$ and $\|a'\|,\|a''\|,\|b'\| \leq L$. Therefore from the above two inequalities we get
      \begin{equation}
      \begin{split}
      &\|[\psi(d_i),a_j]\| \leq 2L\|\psi(d_i) - \psi(c_i')\|+ 2(L+1)\|a_j - a_j'\| + \|[\psi(d_i'),a_j]\|; \\
      &\|[\psi(d_i'),a_j']\| \leq 2(L+1)\|\psi(d_i') - \phi(d_i')\| + \|[\phi(d_i'),a_j']\|.
      \end{split}
      \end{equation}
          Using these approximations we have
\begin{equation}
\begin{split}
              \|[\psi(d_i),a_j]\| &\leq 2L\|\psi(d_i) - \psi(d_i')\| + 2(L+1)\|a_j - a_j'\| + \|[\psi(d_i'),a_j]\| \\
                                &< (4L+2)\delta + \|[\psi(d_i'),a_j]\| \\
                                &\leq (4L+2)\delta + 2(L+1)\|\psi(d_i') - \phi(d_i')\| + \|[\phi(c_i'),a_j']\| \\
                                &< (4L+2)\delta + 2(L+1)\delta + \delta \\
                                &= (6L+5)\delta = \ee.
\end{split}
\end{equation}
  \end{proof}
    \end{prop}

The following will be useful to show that there are many $\cD$-stable embeddings.

\begin{lemma}\label{conjugations}
    Let $\phi: B \simeq B'$ and $\psi: A \simeq A'$ be isomorphisms between unital C*-algebras and let $D$ be a unital, simple, nuclear C*-algebra. Suppose that there is a *-homomorphism $\eta: B' \into A'$ such that there are $(\cF,\ee,\cG)$-embeddings $D \to \eta(B') \subseteq A'$ for all finite subsets $\cF \subseteq D, \cG \subseteq A'$ and $\ee > 0$. Let $\sigma = \psi^{-1} \circ \eta \circ \phi: B \to A$. Then there are $(\cF,\ee,\cG)$-embeddings $D \to \sigma(B) \subseteq A$ for all $\cF \subseteq D, \cG \subseteq A$ finite and $\ee > 0$. 
      \begin{proof}
      The diagram
      \begin{equation}
      \begin{tikzcd}
A \arrow[r, "\psi"]                      & A'                    \\
B \arrow[u, "\sigma"] \arrow[r, "\phi"'] & B' \arrow[u, "\eta"']
\end{tikzcd}
      \end{equation}
      commutes, and so if $\cF \subseteq D, \cG \subseteq A$ are finite, $\ee > 0$ and $\xi: D \to \eta(B') \subseteq A'$ is an $(\cF,\ee,\psi(\cG))$-embedding, then $\psi^{-1} \circ \xi: D \to \psi^{-1}(\eta(B')) \subseteq \psi^{-1}(A') = A$ is an $(\cF,\ee,\cG)$ embedding.  Moreover, from
      \begin{equation}
      \psi^{-1}(\eta(B')) = \psi^{-1}(\eta(\phi(B))) = \sigma(B),
      \end{equation}
      its clear that $\psi^{-1} \circ \xi$ is an $(\cF,\ee,\cG)$ embedding $D \to \sigma(B) \subseteq A$.
      \end{proof}
  \end{lemma}

  \section{Relative intertwinings and $\cD$-stable embeddings}\label{section:dse-embeddings}

  \subsection{Relative intertwinings}\label{ssectoin:dse-relative-intertwinings}  It is well known that a strongly self-absorbing C*-algebra $\cD$ embeds unitally into the central sequence algebra $(\cM(A))_\omega \cap A'$ of a separable C*-algebra $A$ if and only if $A \simeq A \otimes \cD$, where $\cM(A)$ is the multiplier algebra of $A$ (for example, \cite[Theorem 7.2.2(i)]{RordamBook}). We alter the proof to keep track of a subalgebra in order to show that for a unital inclusion $B \subseteq A$ of separable C*-algebras, $\cD \into B_\omega \cap A'$ if and only if there is an isomorphism $\Phi: A \to A \otimes \cD$, which is approximately unitarily equivalent to the first factor embedding, and satisfies $\Phi(B) = B \otimes \cD$. This was initially done for (irreducible) inclusions of $\II_1$ factors in \cite{Bisch90} and commented on in \cite{Izumi04} for $\cD$ being $M_{n^{\infty}},\cO_2,\cO_{\infty}$. The proof we alter is Elliott's intertwining argument, which can be found as a combination of Proposition 2.3.5, Proposition 7.2.1 and Theorem 7.2.2 of \cite{RordamBook}.

  \begin{prop}[Relative intertwining]\label{relativeintertwining}
      Let $A,B,C$ be unital, separable C*-algebras, and let $\phi: A \into C, \theta: B \to A, \psi: B \to C$ be unital *-homomorphisms such that $\phi \circ \theta(B) \subseteq \psi(B)$. Suppose there is a sequence$(u_n) \subseteq \psi(B)_\omega \cap \phi(A)'$ of unitaries such that
        \begin{itemize}
          \item $\dist(v_n^*cv_n,\phi(A)_\omega) \to 0$ for all $c \in C$;
          \item $\dist(v_n^*\psi(b)v_n,\phi \circ \theta(B)_\omega) \to 0$ for all $b \in B$.
        \end{itemize}
        Then $\phi$ is approximately unitarily equivalent to an isomorphism $\Phi: A \simeq C$ such that $\Phi \circ \theta(B) = \psi(B)$.
          \begin{proof}
            Apply the below proposition with $B_m := B, \theta_m := \theta, \psi_m := \psi$ for all $m \in \bN$.
          \end{proof}
    \end{prop}

    \begin{prop}[Countable relative intertwining]\label{countablerelativeintertwining}
        Let $A,B_m,C$ be unital, separable C*-algebras, $m \in \bN$, and $\phi: A \into C,\theta_m: B_m \to A,\psi_m: B_m \to C$ be such that $\phi \circ \theta_m(B_m) \subseteq \psi_m(B_m)$ and $\psi_1(B_1) \subseteq \psi_m(B_m)$. Suppose there is a sequence $(v_n) \subseteq \psi_1(B_1)_\omega \cap \phi(A)'$ of unitaries such that
          \begin{itemize}
            \item $\dist(v_n^*cv_n,\phi(A)_\omega) \to 0$ for all $c \in C$;
            \item $\dist(v_n^*\psi_m(b)v_n,\phi \circ \theta_m(B_m)_\omega) \to 0$ for all $b \in B_m$. 
          \end{itemize}
          Then $\phi$ is approximately unitarily equivalent to an isomorphism $\Phi:A \simeq C$ such that $\Phi \circ \theta_m (B_m) = \psi_m(B_m)$ for all $m \in \bN$. 
            \begin{proof}
              We show that if there are unitaries $(v_n) \subseteq \psi_1(B_1)$ satisfying
                \begin{itemize}
                  \item $[v_n,\phi(a)] \to 0$ for all $a \in A$;
                  \item $\dist(v_n^*cv_n,\phi(A)) \to 0$ for all $c \in C$;
                  \item $\dist(v_n^*\psi_m(b)v_n,\phi \circ \theta_m(B_m)) \to 0$ for all $b \in B_m$,
                \end{itemize}
                then the conclusion holds. Such unitaries can be found using Kirchberg's $\ee$-test (Lemma \ref{lemma:Kirchberg-eptest}).

              Let $(a_n),(b_n^{(m)}),(c_n)$ be dense sequences of $A,B_m,C$ respectively. We can inductively choose $v_n$, forming a subsequence $(v_n)$ of the unitaries above (after reindexing, we are still calling them $v_n$), such that there are $a_{jn} \in A, b_{jn}^{(m)} \in B_m$ with
                \begin{itemize}
                  \item $v_n^*\cdots v_1^*c_jv_1\cdots v_n \approx_{\frac{1}{n}} \phi(a_{jn})$;
                  \item $v_n^*\cdots v_1^*\psi(b_j^{(m)})v_1\cdots v_n \approx_{\frac{1}{n}} \phi \circ \theta_m(b_{jn}^{(m)})$;
                  \item $[v_n,\phi(a_j)] \approx_{\frac{1}{2^n}} 0$;
                  \item $[v_n,\phi(a_{jl})] \approx_{\frac{1}{2^n}} 0$;
                  \item $[v_n,\phi \circ \theta_m(b_j^{(m)})] \approx_{\frac{1}{2^n}} 0$;
                \item $[v_n,\phi \circ \theta_m(b_{jl}^{(m)})] \approx_{\frac{1}{2^n}} 0$,
                \end{itemize}
                where $j,m=1,\dots,n$ and $l=1,\dots,n-1$. Define, for $a \in \{a_n \mid n \in \bN\}$,
\begin{equation}
                \Phi(a) = \lim_n v_1\cdots v_n\phi(a) v_n^*\cdots v_1^*
\end{equation}
                which extends to a *-isomorphism $\Phi: A \simeq C$, as in \cite[Proposition 2.3.5]{RordamBook}. The proof also yields the following useful approximation:
\begin{equation}
                \Phi \circ \theta_m(b_{jn}^{(m)}) \approx_{\frac{1}{2^n}} v_1 \cdots v_n \phi \circ \theta_m(b_{jn}^{(m)})v_n^* \cdots v_1^*
\end{equation}
                for appropriate $n \geq j,m$. 

                We now need to check that $\Phi \circ \theta_m(B_m) = \psi_m(B_m)$. Approximate
\begin{equation}
                 \psi_m(b_j^{(m)}) \approx_{\frac{1}{n}} v_1\cdots v_n \phi \circ \theta_m(b_{jn}^{(m)})v_n^*\cdots v_1^* \approx_{\frac{1}{2^n}} \Phi \circ \theta_m(b_{jn}^{(m)}).
\end{equation}
                This yields $\psi_m(B_m) \subseteq \ov{\Phi \circ \theta_m(B_m)} = \Phi \circ \theta_m(B_m)$. On the other hand for any $\ee > 0$ and $b \in B_m$, we can find $n$ such that
\begin{equation}
                \Phi \circ \theta_m(b) \approx_\ee v_1\cdots v_n \phi \circ \theta_m(b)v_n^*\cdots v_1^* \in \psi_m(B_m)
\end{equation}
                since $v_i \in \psi_1(B_1) \subseteq \psi_m(B_m)$ and $\phi \circ \theta_m(B_m) \subseteq \psi_m(B_m)$. Hence $\Phi \circ \theta_m(B_m) \subseteq \ov{\psi_m(B_m)} = \psi_m(B_m)$. 
            \end{proof}
      \end{prop}

  \subsection{$\cD$-stable embeddings}\label{ssection:dse-d-stable-embeddings}

    \begin{defn}
      Let $\iota: B \into A$ be an embedding and $\cD$ be strongly self-absorbing. We say that $\iota$ is $\cD$-stable (or $\cD$-absorbing) if there exists an isomorphism $\Phi:A \simeq A \otimes \cD$ such that $\Phi \circ \iota(B) = \iota(B) \otimes \cD$. 
    \end{defn}

   We will mostly have interest in the case where $\iota$ corresponds to the inclusion map and $B \subseteq A$ is a subalgebra. In this form, we will say $B \subseteq A$ is $\cD$-stable (or $\cD$-absorbing). Clearly $\iota$ being $\cD$-stable is the same as $\iota(B) \subseteq A$ being $\cD$-stable. We note that we can define the above for any *-homomorphism. Namely, a *-homomorphism $\phi:B \to A$ is $\cD$-stable if $\phi(B) \subseteq A$ is.

   \begin{lemma}\label{trivialinclusion}
       If $\iota: B \into A$ is an embedding, then $\iota \otimes \id_D: B \otimes \cD \into A \otimes \cD$ is $\cD$-stable.
         \begin{proof}
           Let $\phi: D \simeq D \otimes \cD$ be an isomorphism. Then
\begin{equation}
           \Phi := \id_A \otimes \phi: A \otimes \cD \to A \otimes \cD \otimes \cD
\end{equation}
           is an isomorphism with
\begin{equation}
           \Phi(\iota \otimes \id_\cD(B \otimes \cD)) = \left(\iota \otimes \id_\cD(B \otimes \cD)\right) \otimes \cD.
\end{equation}
         \end{proof}
     \end{lemma}

     We note that this is a strengthening of the notion of $\cD$-stability for C*-algebras because if $\iota = \id_A: A \to A$, then $\iota$ is $\cD$-stable if and only if $A$ is $\cD$-stable. This condition is different than the notion of $\cO_2$ or $\cO_{\infty}$-absorbing morphisms discussed in \cite{BGSW22,Gabe20,Gabe19} -- they require sequences from a larger algebra to commute with a smaller algebra, while we require sequences from a smaller algebra to commute with the larger algebra.

   The following adapts  \cite[Theorem 7.2.2]{RordamBook}.
     \begin{theorem}\label{equivDstable}
         Suppose that $B \subseteq A$ is a unital inclusion of separable C*-algebras. If $\cD$ is strongly self-absorbing, then $B \subseteq A$ is $\cD$-stable if and only if there is a unital inclusion $\cD \into B_\omega \cap A'$. 
           \begin{proof}
             Let $\phi: A \to A \otimes \cD$ be the first factor embedding $\phi(a) := a \otimes 1_\cD$. First suppose that $\xi: \cD \into B_\omega \cap A' \simeq (B \otimes 1_\cD)_\omega \cap (A \otimes 1_\cD)'$ is an embedding (so that $\phi(a)\xi(d) \in \phi(A)_\omega$ and $\phi(b)\xi(d) \in \phi(B)_\omega$). Let $\eta: \cD \into (B \otimes \cD)_\omega \cap (A \otimes 1_\cD)'$ be given by $\eta(d) := (1 \otimes d)_n$ and notice that $\xi,\eta$ have commuting ranges. As all endomorphisms of $\cD$ are approximately unitarily equivalent by \cite[Corollary 1.12]{TomsWinter07}, let $(v_n) \subseteq C^*(\xi(\cD),\eta(\cD)) \simeq \cD \otimes \cD$ be such that $v_n^*\eta(d)v_n \to \xi(d)$ for $d \in \cD$. For $b \in B$ and $d \in \cD$, we have
\begin{equation}
\begin{split}
             v_n^*(b \otimes d)v_n &= v_n^*(b \otimes 1_\cD)(1_A \otimes d)v_n^* \\
             &= v_n^*\phi(b)\eta(d)v_n \\
             &= \phi(b)v_n^*\eta(d)v_n \\
             &\to \phi(b)\xi(d) \in \phi(B)_\omega. 
             \end{split}
\end{equation}
             Moreover the same argument shows that, for $a \in A$, we have
\begin{equation}
             v_n^*(a \otimes d)v_n \to \phi(a)\xi(d) \in \phi(A)_\omega.
\end{equation}
             Now $(v_n)$ satisfy the hypothesis of Proposition \ref{relativeintertwining} with $C := A \otimes D$, $\phi$ being the first factor embedding, $\theta: B \to A$ being the inclusion and $\psi: B \simeq B \otimes \cD \subseteq A \otimes \cD = C$ (where this isomorphism exists since if $\cD \into B_\omega \cap A'$, then clearly $\cD \into B_\omega \cap B'$). From this we see that $\phi$ is approximately unitarily equivalent to an isomorphism $\Phi:A \simeq A \otimes \cD$ such that $\Phi(B) = B \otimes \cD$.

             Conversely, if $B \subseteq A$ is $\cD$-stable, let $\Phi: A \simeq A \otimes \cD$ be an isomorphism such that $\Phi(B) = B \otimes \cD$. By \cite[Proposition 1.10(iv)]{TomsWinter07}, we can identify $\cD \simeq \cD^{\otimes \infty}$ and take $\xi: \cD \into B_\omega \cap A'$ to be given by
\begin{equation}
             \xi(d) = (\Phi^{-1}(1_A \otimes 1_\cD^{\otimes n-1} \otimes d \otimes 1_\cD^{\otimes \infty}))_n.
\end{equation}
           \end{proof}
       \end{theorem}

         \begin{cor}
           Let $\iota: B \into A$ be a unital embedding between separable C*-algebras. If $\cD$ is strongly self-absorbing and $\iota$ is $\cD$-stable, then for every intermediate unital C*-algebra $C$ with $\iota(B) \subseteq C \subseteq A$, we have that $\iota(B) \subseteq C$ and $C \subseteq A$ are $\cD$-stable. In particular, $C \simeq C \otimes \cD$ for all such $C$. 
           \begin{proof}
           We have
           \begin{equation}
           \cD \into B_\omega \cap A' \subseteq B_\omega \cap C'
           \end{equation}
           and
           \begin{equation}
           \cD \into B_\omega \cap A' \subseteq C_\omega \cap A'.
           \end{equation}
           \end{proof}
         \end{cor}

         It is not however the case that any isomorphism $\Phi: A \simeq A \otimes \cD$ with $\Phi(B) = B \otimes \cD$ maps $C$ to $C \otimes \cD$.

         \begin{example}\label{noteveryintermediate}
             Let $\cD$ be strongly self-absorbing and consider
\begin{equation}
\begin{split}
                 A &:= \cD \otimes \cD \otimes \cD, \\
                 C_1 &:= \cD \otimes 1_\cD \otimes \cD,\\
                 C_2 &:= 1_\cD \otimes \cD \otimes \cD,\\
                 B &:= 1_\cD \otimes 1_\cD \otimes \cD.
\end{split}
\end{equation}
               If $f: \cD \otimes \cD \to \cD \otimes \cD$ is the tensor flip and $\phi: \cD \simeq \cD \otimes \cD$ is an isomorphism, let
\begin{equation}
                \Phi := f \otimes \phi: A \simeq A \otimes \cD
\end{equation}
               which satisfies $\Phi(B) = B \otimes \cD$ (in particular $B \subseteq A$ is $\cD$-stable). However,
\begin{equation}
               \Phi(C_1) = C_2 \otimes \cD \text{ and  } \Phi(C_2) = C_1 \otimes \cD.
\end{equation}
           \end{example}

           In fact the above example can be generalized to show that for any unital $\cD$-stable inclusion $B \subseteq A$, there is an isomorphism $\Phi: A \simeq A \otimes \cD$ with $\Phi(B) = B \otimes \cD$, and some intermediate algebra $B \subseteq C \subseteq A$ with $\Phi(C) \neq C \otimes \cD$ (obviously we will always have that $\Phi(C) \simeq C \simeq C \otimes \cD$, but equality may not happen).

             \begin{cor}
               Let $B \subseteq A$ be a $\cD$-stable inclusion. There exist a C*-algebra $C$ with $B \subseteq C \subseteq A$ and an isomorphism $\Phi: A \simeq A \otimes \cD$ such that $\Phi(B) = B \otimes \cD$ but $\Phi(C) \neq C \otimes \cD$. 
             \end{cor}

  However, we can always realize $\cD$-stability for countably many intermediate C*-algebras at once using \emph{some} isomorphism $A \simeq A \otimes \cD$.

  \begin{theorem}\label{countablymanysubalgebras}
      Suppose that $B_1 \subseteq B_m \subseteq A$ are unital inclusions of separable C*-algebras (note that we are {\bf not} asking for $(B_m)$ to form a chain). If $\cD$ is strongly self-absorbing and $\cD \into (B_1)_\omega \cap A'$ unitally, there exists an isomorphism $\Phi: A \simeq A \otimes \cD$ such that $\Phi(B_m) = B_m \otimes \cD$ for all $m \in \bN$.
        \begin{proof}
          This is essentially the same proof as Theorem \ref{equivDstable}, except we use the countable relative intertwining (Proposition \ref{countablerelativeintertwining}) in place of Proposition \ref{relativeintertwining}. Let $\xi,\eta$ be as before and let $(v_n) \subseteq C^*(\xi(\cD),\eta(\cD)) \simeq \cD \otimes \cD$ be such that $v_n^*\eta(d)v_n \to \xi(d)$ for $d \in \cD$.
            \begin{itemize}
              \item If $a \in A, d \in \cD, v_n^*(a \otimes d)v_n \to \phi(a)\xi(d) \in \phi(A)_\omega$;
              \item if $b \in B_m, v_n^*(b \otimes d)v_n \to \phi(b)\xi(d) \in \phi(B_m)_\omega$.
            \end{itemize}
            Now with  $\phi: A \to A \otimes \cD$ the first factor embedding,  $\theta_m:B_m \to A$ the inclusion maps, and $\psi_m:B_m \simeq B_m \otimes \cD$ (these exist since $\cD \into (B_1)_\omega \cap A'$ implies that $\cD \into (B_m)_\omega \cap B_m'$), our unitaries satisfy the hypothesis of Proposition \ref{countablerelativeintertwining} and therefore $\phi$ is approximately unitarily equivalent to a *-isomorphism $\Phi: A \simeq A \otimes \cD$ such that $\Phi(B_m) = B_m \otimes \cD$ for all $m$.
        \end{proof}
    \end{theorem}

    The above works since norm ultrapowers have the property that unitaries lift to sequences of unitaries.\footnote{If $u = (u_n) \in A_\omega$, then $\{n \in \bN \mid \|u_n^*u_n - 1\|\|u_nu_n^* - 1\| < 1\} \in \omega$. If $n$ is in the set, replace $u_n$ with the unitary part of its polar decomposition, and replace $u_n$ with 1 otherwise.} Tracial ultrapowers of $\II_1$ von Neumann algebras also have this property.\footnote{The tracial ultrapower of a $\II_1$ von Neumann algebra is again a $\II_1$ von Neumann algebra. Therefore if $u \in \cM^\omega$ is unitary, it is of the form $e^{ia}$ for some $a = a^* \in \cM^\omega$. Lift $a$ to a sequence $(a_n)$ of self-adjoints in $\cM$ and note that $u = (e^{ia_n})$, so that $u$ has a unitary lift.} Consequently if we work with the 2-norm $\|x\|_2 = \tau(x^*x)$ where $\tau$ is the unique trace on a $\II_1$ factor, all of the above arguments with the C*-norm replaced by $\|\cdot\|_2$ will yield back Bisch's result \cite[Theorem 3.1]{Bisch90}, provided we have the appropriate separability conditions. 

    \begin{theorem}
      Let $\cN \subseteq \cM$ be an inclusion of $\II_1$ factors with separable preduals. Then $\cR \into \cN^\omega \cap \cM'$ if and only if there exists an isomorphism $\Phi: \cM \to \cM \ov{\otimes} \cR$ such that $\Phi(\cN) = \cN \ov{\otimes} \cR$.
      \end{theorem}

  \subsection{Existence of $\cD$-stable embeddings}\label{ssection:dse-existence}

  We move to discuss the existence of $\cD$-stable embeddings. First we show that each unital embedding of unital, separable $\cD$-stable C*-algebras is approximately unitarily equivalent to a $\cD$-stable embedding. From this it will follow that there are many $\cD$-stable embeddings. 

  \begin{lemma}\label{adudstable}
      Let $\cD$ be strongly self-absorbing. If $\iota: B \into A$ is a unital, $\cD$-stable inclusion of separable C*-algebras and $u \in U(A)$, then $\Ad_u \circ \iota: B \into A$ is $\cD$-stable.
        \begin{proof}
          Apply Lemma \ref{aduaapproximateembeddings}. 
        \end{proof}
    \end{lemma}

    \begin{prop}\label{allembeddingsareue}
        Let $\cD$ be strongly self-absorbing, $A,B$ be unital separable $\cD$-stable C*-algebras and let $\iota: B \into A$ be an embedding. Then $\iota$ is approximately unitarily equivalent to a $\cD$-stable embedding $B \into A$.
          \begin{proof}
            As $A,B$ are $\cD$-stable, there are isomorphisms
\begin{equation}
             \phi: B \simeq B \otimes \cD \text{ and }\psi: A \simeq A \otimes \cD
\end{equation}
            which are approximately unitarily equivalent to the first factor embeddings $b \mapsto b \otimes 1_\cD, b \in B$ and $a \mapsto a \otimes 1_\cD, a \in A$ respectively. As $\iota \otimes \id_\cD: B \otimes \cD \into A \otimes \cD$ is $\cD$-stable by Lemma \ref{trivialinclusion},
\begin{equation}
            \sigma := \psi^{-1} \circ (\iota \otimes \id_\cD) \circ \phi: B \into A
\end{equation}
            is $\cD$-stable by Lemma \ref{conjugations}. Now we show that $\sigma$ is approximately unitarily equivalent to $\iota$. Let $\cF \subseteq B$ be finite and $\ee > 0$. Let $u \in U(B \otimes \cD)$ be such that $u^*(b \otimes 1_\cD)u \approx_\frac{\ee}{2} \phi(b)$ for $b \in \cF$ and $v \in U(A \otimes \cD)$ be such that $v^*(\iota(b) \otimes 1_\cD)v \approx_{\frac{\ee}{2}} \psi \circ \iota (b)$ for $b \in \cF$. Set $w = \psi^{-1}(\iota \otimes \id_\cD(u))^*\psi^{-1}(v) \in U(A)$. Then for $b \in \cF$,
\begin{equation}
\begin{split}
                w^*\sigma(b)w &= \psi^{-1}(v)^*\psi^{-1}(\iota \otimes \id_\cD(u\phi(b)u^*))\psi^{-1}(v) \\
                              &\approx_{\frac{\ee}{2}} \psi^{-1}(v)^*\psi^{-1}(\iota \otimes \id_\cD(b \otimes 1_\cD))\psi^{-1}(v) \\
                              &= \psi^{-1}(v)^*\psi^{-1}(\iota(b) \otimes 1_\cD)\psi^{-1}(v) \\
                              &\approx_{\frac{\ee}{2}} \psi^{-1}(\psi(\iota(b))) \\
                              &= \iota(b).
\end{split}
\end{equation}

          \end{proof}
      \end{prop}

      \begin{cor}\label{pointnormdensity}
          Let $\cD$ be strongly self-absorbing. The set of $\cD$-stable embeddings $B \into A$ of unital, separable, $\cD$-stable C*-algebras is point-norm dense in the set of embeddings $B \into A$.
            \begin{proof}
              Every embedding is approximately unitarily equivalent to a $\cD$-stable embedding. As $\cD$-stability of an embedding is preserved if one composes with $\Ad_u$, it follows that every embedding is the point-norm limit of $\cD$-stable embeddings.
            \end{proof}
        \end{cor}

        \begin{remark}\label{homomorphisminsteadofembedding}
          We note that it is not actually necessary that $\iota$ is an embedding. If $\pi: B \to A$ is any unital *-homomorphism between unital, separable, $\cD$-stable C*-algebras, then $\pi$ is approximately unitarily equivalent to a *-homomorphism $\pi': B \to A$ such that $\pi'(B) \subseteq A$ is $\cD$-stable. Consequently the set of unital *-homomorphisms $\pi: B \to A$ with $\pi(B) \subseteq A$ being $\cD$-stable is in fact dense in the set of unital *-homomorphisms $B \to A$.
        \end{remark}

        Later on, there will be some examples of non-$\cD$-stable embeddings between $\cD$-stable C*-algebras. Consequently, despite the fact $\cD$-stable embeddings are point-norm dense, the set of $\cD$-stable embeddings need not coincide with the set of all embeddings $B \into A$. Another clear consequence is that despite $\cD$-stability of an embedding being closed under conjugation by a unitary, it is not true that it is preserved under approximate unitary equivalence (in fact, the examples in question show that $\cD$-stability is not even preserved under asymptotic unitary equivalence). We finish with a corollary about embeddings into the Cuntz algebra  $\cO_2$ \cite{Cuntz77}. 

          \begin{cor}
            Let $B$ be a unital, separable, exact $\cD$-stable C*-algebra, where $\cD$ is strongly self-absorbing. Then there is a $\cD$-stable embedding $B \into \cO_2$. 
              \begin{proof}
                As $\cD$ is unital, simple, separable and nuclear by \cite[Section 1.6]{TomsWinter07}, $\cO_2 \simeq \cO_2 \otimes \cD$ and $B \into \cO_2$ unitally by Theorem 3.7 and Theorem 2.8 of \cite{KirchbergPhillips00} respectively. The above results then yield a $\cD$-stable embedding $B \into \cO_2$. 
              \end{proof}
          \end{cor}

  We include this last result about the classification of morphisms via functors. 

  \begin{theorem}
      Let $\cD$ be strongly self-absorbing and let $F$ be a functor from a class of unital, separable, $\cD$-stable C*-algebras satisfying the following.
        \begin{enumerate}
          \item[(E)] If there exists a morphism $\Phi: F(B) \to F(A)$, then there exists a *-homomorphism $\phi: B \to A$ such that $F(\phi) = \Phi$.
          \item[(U)] If $\phi,\psi: B \to A$ are *-homomorphisms which are approximately unitarily equivalent, then
\begin{equation}
             F(\phi) = F(\psi).
\end{equation}
        \end{enumerate}
        Then whenever there is a morphism $\Phi: F(B) \to F(A)$, there exists $\phi: B \to A$ such that  $F(\phi) = \Phi$ and $\phi(B)  \subseteq A$ is $\cD$-stable.
        Moreover, $\phi$ is unique up to approximate unitary equivalence.
          \begin{proof}
            By the existence (E), there exists a *-homomorphism $\phi: B \to A$. Now by Proposition \ref{allembeddingsareue} (Remark \ref{homomorphisminsteadofembedding} allows us to work with general *-homomorphisms), there exists a *-homomorphism $\phi': B \to A$ which is approx\hyp{}imately unitarily equivalent to $\phi$ and $\phi'(B) \subseteq A$ is $\cD$-stable.
            Uniqueness (U) gives that this is unique up to approximate unitary equivalence. 
          \end{proof}
    \end{theorem}

  \subsection{Permanence properties}\label{ssection:dse-permanence}
  We now discuss some permanence properties.

    \begin{lemma}
Let $\cD$ be strongly self-absorbing.
Suppose that $\iota_i: B_i \into A_i, i=1,2$ are $\cD$-stable inclusions. Then $\iota_1 \oplus \iota_2: B_1 \oplus B_2 \into A_1 \oplus A_2$ is $\cD$-stable.
        \begin{proof}
          Let $\Phi_i: A_i \simeq A_i \otimes \cD$ be isomorphisms such that $\Phi_i\circ \iota_i(B_i) = \iota_i(B_i) \otimes \cD$ and consider
\begin{equation}
          \Phi: A_1 \oplus A_2 \simeq (A_1 \oplus A_2) \otimes \cD
\end{equation}
          given by the composition
\begin{equation}
          \begin{tikzcd}
A_1 \oplus A_2 \arrow[r, "\Phi_1 \oplus \Phi_2"] & (A_1 \otimes \cD) \oplus (A_2 \otimes \cD) \arrow[r, "\simeq"] & (A_1 \oplus A_2) \otimes \cD
\end{tikzcd}
\end{equation}
where the last isomorphism follows from (finite) distributivity of the min-tensor. Then we see that
\begin{equation}
\Phi(\iota_1(B_1) \oplus \iota_2(B_2)) = \left(\iota_1(B_1) \oplus \iota_2(B_2)\right) \otimes \cD.
\end{equation}
        \end{proof}
    \end{lemma}

      \begin{lemma}
Let $\cD$ be strongly self-absorbing.        
        Suppose that $\iota_i: B_i \into A_i, i=1,2$ are inclusions and that at least one of $\iota_1$ or $\iota_2$ is $\cD$-stable. Then $\iota_1 \otimes \iota_2: B_1 \otimes B_2 \into A_1 \otimes A_2$ is $\cD$-stable. 
          \begin{proof}
            We prove this if $\iota_2$ is $\cD$-stable, and a symmetric argument will yield the result if $\iota_1$ is. Let $\Phi_2: A_2 \simeq A_2 \otimes \cD$ be such that $\Phi_2 \circ \iota_2(B_2) = \iota(B_2) \otimes \cD$. Taking
\begin{equation}
            \Phi := \id_{A_1} \otimes \Phi_2: A_1 \otimes A_2 \simeq A_1 \otimes A_2 \otimes \cD,
\end{equation}
            we have that
\begin{equation}
            \Phi(\iota_1(B_1) \otimes \iota_2(B_2)) = \iota_1(B_1) \otimes \iota_2(B_2) \otimes \cD.
\end{equation}
          \end{proof}
      \end{lemma}
      
        \begin{prop}
          Let $\cD$ be strongly self-absorbing.
Suppose that we have increasing sequences of unital separable C*-algebras $(B_n)$ and $(A_n)$ such that $B_n \subseteq A_n$ unitally.
Let $B = \ov{\cup_nB_n}$ and $A = \ov{\cup_n A_n}$. If $B_n \subseteq A_n$ is $\cD$-stable for all $n$, then $B \subseteq A$ is $\cD$-stable.
  \begin{proof}
    This follows from Proposition \ref{approximateembeddinginductivelimit}, together with Lemma \ref{lem:local-equiv-global} and Theorem \ref{equivDstable}.
  \end{proof}
        \end{prop}

        Lastly we'll discuss unital inclusions $B \subseteq A$ of $C(X)$ algebras, where $X$ is a compact Hausdorff space. We show that if $X$ has finite covering dimension, then such an inclusion is $\cD$-stable if and only if the inclusion $B_x \subseteq A_x$ along each fibres is $\cD$-stable. 

        \begin{lemma}\label{imageofinclusion}
            Let $\cD$ be strongly self-absorbing.
Suppose that $B_i \subseteq A_i$ are unital inclusions, for $i=1,2$, and $\psi: A_1 \to A_2$ is a surjective *-homomorphism such that $\psi(B_1) = B_2$. If $B_1 \subseteq A_1$ is $\cD$-stable, then so is $B_2 \subseteq A_2$.
              \begin{proof}
                We note that $\psi$ induces a *-homomorphism
\begin{equation}
                \tilde{\psi}: (B_1)_\omega \cap A_1' \to (B_2)_\omega \cap A_2'
\end{equation}
                and consequently if $\xi: \cD \into (B_1)_\omega \cap A_1'$, we have a unital *-homomorphism
\begin{equation}
                \eta := \tilde{\psi} \circ \xi: \cD \to (B_2)_\omega \cap A_2'.
\end{equation}
                $\eta$ is automatically injective since $\cD$ is simple. 
              \end{proof}
          \end{lemma}

Rephrasing the above in terms of commutative diagrams, it says that if we have a commutative diagram
\begin{equation}
\begin{tikzcd}
A_1 \arrow[r, two heads]                 & A_2                 \\
B_1 \arrow[r, two heads] \arrow[u, hook] & B_2 \arrow[u, hook]
\end{tikzcd}
\end{equation}
where the left inclusion is $\cD$-stable, then the right inclusion is $\cD$-stable as well. 

Now we consider many of the results discussed in \cite[Section 4]{HirshbergRordamWinter07}, except for inclusions of C*-algebras.

\begin{defn}
      Let $X$ be a compact Hausdorff space. A $C(X)$-algebra is a C*-algebra $A$ endowed with a unital *-homomorphism $C(X) \to \cZ(\cM(A))$, where $\cZ(\cM(A)$ is the center of the multiplier algebra $\cM(A)$ of $A$. 
    \end{defn}

    If $Y \subseteq X$ is a closed subset, we set $I_Y := C_0(X \sm Y)A$, which is a closed two-sided ideal in $A$. We denote $A_Y := A/A_Y$ and the quotient map  $A \to A_Y$ by $\pi_Y$. For an element $a \in A$, we write $a_Y := \pi_Y(a)$ and if $Y$ consists of a single point $x$, we write $A_x,I_x,\pi_x$ and $a_x$. We say that $A_x$ is the fibre of $A$ at $x$. We note that $A_X = A$.
    
    If $B \subseteq A$ is a unital inclusion and $\theta:_A: C(X) \to A, \theta_B: C(X) \to B$ are morphisms witness $A$ and $B$ as $C(X)$-algebras, respectively, we say that $B \subseteq A$ is an inclusion of $C(X)$-algebras if
    \begin{equation}
    \begin{tikzcd}
B \arrow[r, hook]                                  & A \\
C(X) \arrow[u, "\theta_B"] \arrow[ru, "\theta_A"'] &  
\end{tikzcd}
    \end{equation}
    commutes. Note that $\theta_B(B) \subseteq \cZ(A)$. We note that when discussion an inclusion of fibres $B_Y \subseteq A_Y$ we are considering $B_Y := \pi^A_Y(B) \subseteq \pi^A_Y(A) =: A_Y$, where $\pi_Y^A: A \to A_Y$ is the associated quotient map.

  \begin{remark}[Upper semi-continuity]\label{rem:up-semi-cts}
  In \cite[Section 1.3]{HirshbergRordamWinter07}, it was pointed out that the norm on a $C(X)$-algebra $A$ is upper semi-continuous. This meaning that, fixing some $a \in A$, the function $x \mapsto \|a_x\|$ from $X$ to $\bR$ is upper semi-continuous (as it is the infimum of a family of continuous functions), and consequently the set $\{x \in X \mid \|a_x\| < \ee\} \subseteq X$ is open for all $a \in A$ and $\ee > 0$. 
  \end{remark}

We note that Lemma \ref{imageofinclusion} gives that if $B\subseteq A$ is $\cD$-stable and $Y \subseteq X$ is closed, then $B_Y \subseteq A_Y$ is automatically $\cD$-stable as well since we have the commuting diagram
\begin{equation}
\begin{tikzcd}
A \arrow[r, "\pi_Y"]                     & A_Y                \\
B \arrow[u, hook] \arrow[r, "\pi_Y|_B"'] & B_Y. \arrow[u, hook]
\end{tikzcd}
\end{equation}
The converse needs a bit of work. This is the embedding-analogue of the beginning of \cite[Section 4]{HirshbergRordamWinter07}. We discuss how the proofs can be adapted and often omit approximations that were otherwise done there. We want a version of \cite[Lemma 4.5]{HirshbergRordamWinter07}, which is a result about \emph{gluing} c.c.p. maps together along fibres. In our setting, we are only interested in u.c.p. maps, and we want to show that if we \emph{glue} two u.c.p. maps together whose images are contained in some $C(X)$-subaglebra $B$, then the \emph{glued} map also has image contained in $B$. We borrow their Definition 4.2.

\begin{defn}
Let $A$ be a unital $C(X)$-algebra, for a compact Hausdorff space $X$, and let $D$ be a unital C*-algebra. Let $\phi: D \to A$ is a u.c.p. map and $Y \subseteq X$ a closed subset. If $\cF \subseteq D, \cG \subseteq A$ are finite and $\ee > 0$, we say that $\phi$ is $(\cF,\ee,\cG)$-good for $Y$ if
\begin{enumerate}
\item $([\phi(d),a])_Y \approx_\ee 0$ and
\item $\phi(dd')_Y \approx_\ee \phi(d)_Y\phi(d')_Y$
\end{enumerate}
whenever $d,d' \in \cF$ and $a \in \cG$. If $X = [0,1]$, $Y \subseteq X$ is a closed interval, $\cF' \supseteq \cF$ is another finite set and $0 < \ee' < \ee$, we say that $\phi$ is $(\cF,\ee,\cG;\cF',\ee')$-good for $Y$ if $\phi$ is $(\cF,\ee,\cG)$-good for $Y$ and there exists some closed neighbourhood $V$ of the endpoints of $Y$ such that $\phi$ is $(\cF',\ee',\cG)$-good for $V$. 
\end{defn}

First we need a lemma that follows as a consequnce of $\cD$-stability. It is the embedding analogue of \cite[Proposition 4.1]{HirshbergRordamWinter07}.

\begin{lemma}\label{lem:kappa-mu}
Let $\cD$ be strongly self-absorbing, and $B \subseteq A$ be a unital, $\cD$-stable inclusion of separable C*-algebras. Then for any $\cG \subseteq A$ finite and $\ee > 0$, there exist unital *-homomorphisms $\kappa: A \to A$ and $\mu: \cD \to B$ such that
\begin{enumerate}
\item $\kappa(B) \subseteq B$,
\item $[\kappa(A),\mu(\cD)] = 0$,
\item $\kappa(a) \approx_\ee a$ for all $a \in \cG$.
\end{enumerate}
\begin{proof}
The proof is essentially the same as the proof of (a) $\Rightarrow (c)$ in \cite[Proposition 4.1]{HirshbergRordamWinter07}. As $B \subseteq A$ is $\cD$-stable, let us identify $B \subseteq A$ with $B \otimes \cD \subseteq A \otimes \cD$. As $\cD$ is strongly self-absorbing,  \cite[Theorem 2.3]{TomsWinter07} gives a sequence $(\phi_n)$ of *-homomorphisms $\phi_n: \cD \otimes \cD \to \cD$ such that
\begin{equation}
\phi_n(d \otimes 1_\cD) \to d \text{ for all } d \in \cD.
\end{equation}
Define $\kappa_n: A \otimes \cD \to A \otimes \cD$ by
\begin{equation}
\kappa_n:= (\id_A \otimes \phi) \circ (\id_A \otimes \id_\cD \otimes 1_\cD), 
\end{equation}
and $\mu_n: \cD \to B \otimes \cD$ by
\begin{equation}
\mu_n := (\id_B \otimes \phi_n) \circ (1_A \otimes 1_\cD \otimes \id_\cD).
\end{equation}
Then taking $n$ large enough and letting $\kappa$ and $\mu$ be $\kappa_n$ and $\mu_n$ respectively, its clear that $\kappa(B \otimes \cD) \subseteq B \otimes \cD$, $[\kappa(A),\mu(\cD)] = 0$ and that $\kappa(a) \approx_\ee a$ whenever $a$ is in some prescribed finite subset $\cG \subseteq A$ and $\ee > 0$ is some prescribed error. 
\end{proof}
\end{lemma}

\begin{lemma}
Let $\cD$ be strongly self-absorbing and $A$ be a unital, separable $C([0,1])$-algebra. Suppose $\cF \subseteq \cD, \cG \subseteq A$ are finite self-adjoint subsets of contractions with $1_\cD \in \cF$. Suppose that we have points $0 \leq r < s < t \leq 1$ and two u.c.p. maps $\rho,\sigma: \cD \to A$ which are $(\cF,\ee,\cG)$-good for $[r,s],[s,t]$ respectively. Suppose that $A_s$ is $\cD$-stable. 

Then there are u.c.p. maps $\rho',\sigma': \cD \to A$ which are $(\cF,\ee,\cG)$-good for $[r,s],[s,t]$ respectively, and u.c.p. maps $\nu_{\rho'},\nu_{\sigma'}: \cD \to A,\mu_{\rho'},\mu_{\sigma'}: \cD \otimes \cD \to A$ such that $\nu_{\rho'},\nu_{\sigma'}$ are $(\cF,3\ee,\cG)$-good for some interval $I \subseteq (r,t)$ containing $s$ in its interior, and such that for any $a \in \cG, d,d' \in \cF$, we have
\begin{enumerate}
\item $([\rho'(d),\nu_{\rho'}(d')])_I \approx_{2\ee} 0$
\item $([\sigma'(d),\nu_{\sigma'}(d')])_I \approx_{2\ee} 0$
\item $\rho'(d)_I\nu_{\rho'}(d')_I \approx_\ee \mu_{\rho'}(d \otimes d')_I$
\item $\sigma'(d)_I\nu_{\sigma'}(d')_I \approx_\ee \mu_{\sigma'}(d \otimes d')_I$
\item $\nu_{\rho'}(d)_I \approx_{2\ee} \nu_{\sigma'}(d)_I$.
\end{enumerate} 
If $\rho,\sigma$ are $(\cF,\ee,\cG;\cF',\ee)$-good for $[r,s],[s,t]$ respectively, for some finite $\cF' \supseteq \cF$ set of contractions and for some $0 < \ee' < \ee$, then we can arrange so that $\rho',\sigma',\nu_{\rho'},\nu_{\sigma'}$ are $(\cF',3\ee',\cG)$-good for the interval $I$, and that the above five conditions hold with $\ee'$ in place of $\ee$ and $\cF'$ in place of $\cF$. 

Moreover, if $B \subseteq A$ is a unital inclusion of $C([0,1])$-algebras such that $\rho(\cD) \subseteq B, \sigma(\cD) \subseteq B$ and $B_s \subseteq A_s$ is $\cD$-stable, then the images of all $\rho',\sigma',\mu_{\rho'},\mu_{\sigma'}$ are contained in $B$ (as are the images of $\nu_{\rho'}$ and $\nu_{\sigma'}$).
\begin{proof}
This is \cite[Lemma 4.4]{HirshbergRordamWinter07}, except we've replaced c.c.p. maps with u.c.p. maps. One can easily check that the resulting maps are u.c.p. maps. 

As for the ``moreover'' part, which is the only addition besides the unitality, we outline the definitions of these maps to show that the images of $\rho',\sigma',\mu_{\rho'},\mu_{\sigma'}$ are contained in $B$. As $B_s \subseteq A_s$ is $\cD$-stable, we can find $\kappa: A_s \to A_s$ and $\mu: \cD \to B_s$ as in Lemma \ref{lem:kappa-mu}, where $\kappa(a_s) \approx a_s$ for an appropriate error whenever $a \in \cG$. We use Choi-Effros to find u.c.p. lifts $\tilde{\rho},\tilde{\sigma}: \cD \to B$ for the maps $\kappa \circ \pi_s \circ \rho$ and $\kappa \circ \pi_s \circ \sigma$ respectively (note that $\kappa \circ \pi_s \circ \rho$ and $\kappa \circ \pi_s \circ \sigma$ lie in $B_s$, which is a *-homomorphism image of $B$). One then defines piece-wise linear functions $f,g: [0,1] \to [0,1]$ which attain both values 0 and 1 at the end points (their definition is not important to show the ``moreover'' part). Then $\rho',\sigma'$ are defined as
\begin{equation}
\rho'(d) := (1-f)\cdot \rho(d)  + f\cdot\tilde{\rho}(d) \text{ and } \sigma'(d) := (1-g)\cdot\sigma(d) + g\cdot \tilde{\sigma}(d)
\end{equation}
Clearly $\rho',\sigma'$ take values in $B$ as $\rho,\tilde{\rho},\sigma,\tilde{\sigma}$ all do and $(1-f),f,(1-g),g$ are in $B$. Now we define u.c.p. maps $\tilde{\mu}_{\rho'},\tilde{\mu}_{\sigma'}: \cD \otimes \cD \to B_s$ by
\begin{equation}
\tilde{\mu}_{\rho'}(d \otimes d') := \rho'(d)_s\mu(d') \text{ and }\tilde{\mu}_{\sigma'}(d \otimes d') := \sigma'(d)_s\mu(d').
\end{equation}
Now by Choi-Effros, we can take u.c.p. lifts $\mu_{\rho'}$ and $\mu_{\sigma'}$ of $\tilde{\mu}_{\rho'}$ and $\mu_{\sigma'}$, respectively. As the images of $\tilde{\mu}_{\rho'}$ and $\mu_{\sigma'}$ lie in $B_s$, the images of $\mu_{\rho'}$ and $\mu_{\sigma'}$ will lie in $B$. 
\end{proof}
\end{lemma}

\begin{lemma}\label{lem:exist-F-ee}
Let $A$ be a unital, separable $C([0,1])$-algebra. Suppose $\cF \subseteq \cD, \cG \subseteq A$ are finite self-adjoint subsets with $1_\cD \in \cF$ and $\ee > 0$. There exists $0 < \ee' < \ee$ and a finite subset $\cF' \supseteq \cF$ such that if $\rho,\sigma: \cD \to A$ are u.p.c. maps and $0 \leq r < s < t \leq 1$ are points such that $\rho$ is $(\cF,\ee,\cG;\cF',\ee')$-good for $[r,s]$,  $\sigma$ is $(\cF,\ee,\cG;\cF',\ee')$-good for $[s,t]$ and $A_s$ is $\cD$-stable, then there is a u.c.p. map $\psi: \cD \to A$ which is $(\cF,\ee,\cG;\cF',\ee')$-good for $[r,t]$. 

Moreover, if $B \subseteq A$ is a unital inclusion of $C([0,1])$-algebras such that $\rho(\cD) \subseteq B$, $\sigma(\cD) \subseteq B$ and $B_s \subseteq A_s$ is $\cD$-stable, then $\psi(\cD) \subseteq B$.
\begin{proof}
The first part is \cite[Lemma 4.5]{HirshbergRordamWinter07}, except we've replaced c.c.p. maps with u.c.p. maps. One has to check that the resulting $\psi$ is unital, but this follows easily if $\rho$ and $\sigma$ are.

We outline the construction of $\psi$ to show unitality, as it will also be useful to show the ``moreover'' part, which is the only real addition. Let $u \in C([0,1],\cD \otimes \cD)$ be a path of unitaries such that $u_0 = 1_{\cD \otimes \cD}$ and
\begin{equation}
u_1(d \otimes 1_\cD)u_1^* \approx_{\frac{\ee}{4}} 1_\cD \otimes d. 
\end{equation}
We replace $\rho,\sigma$ with $\rho',\sigma'$ as in the above lemma and this yields u.c.p. maps $\mu_\rho,\mu_\sigma$ satisfying the hypotheses above for some interval $I \subseteq (s,t)$ with $s$ in its interior. Define
\begin{equation}
\phi_\rho,\phi_\sigma: C([0,1]) \otimes \cD \otimes \cD \to A
\end{equation}
by
\begin{equation}
\begin{split}
\phi_\rho(f \otimes d \otimes d') &:= f\cdot \mu_\rho(d \otimes d') \\
\phi_\sigma(f \otimes d \otimes d') &:= f\cdot \mu_\sigma(d \otimes d'). 
\end{split}
\end{equation}
Note that these maps are unital. Take piece-wise linear functions $h_1,h_2,h_3,h_4: [0,1] \to [0,1]$ which sum to 1 (their specific form does not matter to show unitality of $\psi$ nor the ``moreover'' part) and  $g_\rho,g_\sigma: [0,1] \to [0,1]$ which sum to 1 (again, their specific form does not matter to show unitality of $\psi$ nor the ``moreover'' part). Define unitaries $u_\rho,u_\sigma \in C([0,1]) \otimes \cD \otimes \cD \simeq C([0,1],\cD \otimes \cD)$ by 
\begin{equation}
u_{\rho x} := u_{g_\rho(x)} \text{ and } u_{\sigma x} := u_{g_\sigma(x)}.
\end{equation}
Now define $\zeta_\rho,\zeta_\sigma: \cD \to A$ by
\begin{equation}
\begin{split}
\zeta_\rho(d) &:=  \phi_\rho(u_\rho(1_{C([0,1])} \otimes d \otimes 1_\cD)u_\rho^*) \\
\zeta_\sigma(d) &:= \phi_\sigma(u_\sigma(1_{C([0,1])} \otimes d \otimes 1_\cD)u_\sigma^*),
\end{split}
\end{equation}
which are clearly unital. Finally the map $\psi: \cD \to A$ is defined by
\begin{equation}
\psi(d) := h_1\cdot\rho(d) + h_2\cdot \zeta_\rho(d) + h_3\cdot \zeta_\sigma(d) + h_4\cdot \sigma(d).
\end{equation}
Clearly $\psi$ is unital. 

Now for the ``moreover'' part. If $\rho(\cD) \subseteq B$ and $\sigma(\cD) \subseteq B$, clearly the first and fourth terms in the definition of $\psi$ will lie in $B$. So it suffices to show that $\zeta_\rho(\cD) \subseteq B$ and $\zeta_\sigma(\cD) \subseteq B$, and for this it suffices to show that $\mu_\rho(\cD \otimes \cD) \subseteq B$ and $\mu_\sigma(\cD \otimes \cD) \subseteq B$ (since $h_1,h_2,h_3,h_4$ all lie in $B$). But this follows from the ``moreover'' part of the previous lemma. 
\end{proof}
\end{lemma}

  With this, we get the analogue of their Theorem 4.6, the proof being essen\hyp{}tially the same as well, except we insist that the our u.c.p. maps commute with a prescribed finite subset of $A$. 
  
\begin{prop}
Let $\cD$ be strongly self-absorbing, and $X$ be a compact Hausdorff space with finite covering dimension. Suppose that $B \subseteq A$ is a unital inclusion of $C(X)$-algebras. Then $B_x \subseteq A_x$ is $\cD$-stable for all $x \in X$ if and only if $B \subseteq A$ is $\cD$-stable. 
          \begin{proof}
            As previously mentioned, if $B \subseteq A$ is $\cD$-stable, then $B_x \subseteq A_x$ is $\cD$-stable for all $x$.

            For the converse, the proof is essentially the same as \cite[Theorem 4.6]{HirshbergRordamWinter07}. Using the arguments there, one can simplify to the case where we can argue this for $C([0,1])$-algebras (by using \cite[Theorem V.3]{HurewiczWallman}, which says that a compact space of dimension $\leq n$ is homeomorphic to a subset of $[0,1]^{2n+1}$, and then working component-wise). Now for $\cF \subseteq \cD, \cG \subseteq A$ and  $\ee > 0$, let $\cG_x := \{a_x \mid a \in \cG\}$. Without loss of generality suppose that $\cF^*=\cF,\cG^*=\cG$ and that $1_\cD \in \cF$. Let $\cF',\ee'$ be as in Lemma \ref{lem:exist-F-ee}. 
            
            By $\cD$-stability of the inclusion $B_x \subseteq A_x$ there are u.c.p. $(\cF',\ee',\cG_x)$-embed\hyp{}dings $\psi_x: \cD \to B_x \subseteq A_x$ which lift by Choi-Effros to u.c.p. maps $\psi_x': \cD \to B$. The norm is upper semi-continuous (Remark \ref{rem:up-semi-cts}), and this yields intervals $I_x \subseteq [0,1]$ such that $\psi_x'$ is $(\cF',\ee',\cG)$-good for $\ov{I_x}$.  Note that $\psi_x'$ being $(\cF',\ee',\cG)$-good for the whole of $I_x$ implies that it is $(\cF,\ee,\cG;\cF',\ee')$-good for $\ov{I_x}$.  
             Compactness then allows us to split the interval as
\begin{equation}
            0 = t_0 < t_1 < \dots < t_n = 1
\end{equation}
            and to take $\psi_i: \cD \to B$ u.c.p.  which are $(\cF,\ee,\cG;\cF',\ee')$-good for $[t_{i-1},t_i]$ for $i=1,\dots,n$ ($\psi_i = \psi_x'$ for some $x \in [0,1]$). Now by repeatedly using the gluing lemma (Lemma \ref{lem:exist-F-ee}) to glue these maps together, we can find a u.c.p. map $\psi: \cD \to B$ which is an $(\cF,\ee,\cG)$-embedding. 
          \end{proof}
\end{prop}

  \section{Crossed products}\label{section:dse-crossed-products}

  In this section we discuss how inclusions coming from non-commutative dynamics fit into the framework of a tensorially absorbing inclusions. We'll shortly discuss group actions $G \act^\alpha A$ with Rokhlin properties and consider the inclusion of a C*-algebra in its crossed product $A \subseteq A \rtimes_\alpha G$, as well as the inclusion of the fixed point subalgebra of the action in the C*-algebra $A^\alpha \subseteq A$. We then discuss diagonal inclusions associated to certain group actions. 

  This first result says that if we have an isomorphism $A \simeq A \otimes \cD$ which is $G$-equivariant with respect to an action point-wise fixing the right tensor factor, then the corresponding inclusion $A \subseteq A \rtimes_{r,\alpha} G$ is $\cD$-stable.  

  \begin{prop}
    Let $G \act^\alpha A$ be an action of a countable discrete group on a unital separable C*-algebra. Suppose that $\alpha \simeq \alpha \otimes \id_\cD$, that is, there is an isomorphism $\Phi: A \simeq A \otimes \cD$ such that
\begin{equation}
 \begin{tikzcd}
A \arrow[d, "\alpha_g"'] \arrow[r, "\Phi"] & A \otimes \cD \arrow[d, "\alpha_g \otimes \id_\cD"] \\
A \arrow[r, "\Phi"']                       & A \otimes \cD                                      
\end{tikzcd}
\end{equation}
commutes for all $g \in G$. Then $A \subseteq A \rtimes_{r,\alpha} G$ is $\cD$-stable. 
      \begin{proof}
          Let $\psi: \cD \simeq \cD^{\otimes \infty}$ and let $\phi_n: \cD \to \cD^{\otimes \infty}$ be the $n$th factor embedding:
\begin{equation}
            \phi_n(d)  := 1_\cD^{\otimes n-1} \otimes d \otimes 1_\cD^{\otimes \infty}.
\end{equation}
            We claim that $\xi(d) := (\Phi^{-1}(1_A \otimes \psi^{-1} \circ \phi_n(d)))_n: \cD \to A_\omega$ is an embedding such that $\xi(\cD) \subseteq A_\omega \cap A'$ and $(\alpha_g)_\omega \circ \xi = \xi$ for all $g \in G$ -- that is, $\xi$ is an embedding $\cD \into A_\omega \cap (A \rtimes_{r,\alpha} G)'$. The first claim is obvious, so we prove the second. We have
\begin{equation}
\begin{split}
                \|\alpha_g(\Phi^{-1}(1_A &\otimes \psi^{-1}(\phi_n(d)))) -  \Phi^{-1}(1_A \otimes \psi^{-1}(\phi_n(d)))\| \\
                                         &= \|\Phi \circ \alpha_g(\Phi^{-1}(1_A \otimes \psi^{-1}(\phi_n(d)))) -  \Phi(\Phi^{-1}(1_A \otimes \psi^{-1}(\phi_n(d))))\| \\
                                         &= \|\alpha_g \otimes \id_\cD(1_A \otimes \psi^{-1}(\phi_n(d))) - 1_A \otimes \psi^{-1}(\phi_n(d)))\| \\
                                         &=0. 
\end{split}
\end{equation}
      \end{proof}
  \end{prop}

  The next lemma of note is the following. 

    \begin{lemma}\label{lem:fixed-point-in-crossed-product-D-stable}
      Suppose that $G \act^\alpha A$ is an action of a finite group on a unital separable C*-algebra $A$ such that $A \subseteq A \rtimes_\alpha G$ is $\cD$-stable. Then $A^\alpha \subseteq A \rtimes_\alpha G$ is $\cD$-stable. In particular, if $A \subseteq A\rtimes_\alpha G$ is $\cD$-stable, then $C \simeq C \otimes \cD$ whenever $A^\alpha \subseteq C \subseteq A \rtimes_\alpha \cD$.
        \begin{proof}
        For an element $(x_n) \in A_\omega \cap (A \rtimes_\alpha G)'$, an easy averaging argument shows that
\begin{equation}
        (x_n) = \left(\frac{1}{|G|}\sum_{g \in G} \alpha_g(x_n)\right)
\end{equation}
in $A_\omega$, and the right is clearly point-wise fixed by $\alpha_g$ for all $g \in G$. So $A_\omega \cap (A \rtimes_\alpha G)'$ is actually equal to $(A^\alpha)_\omega \cap (A \rtimes_\alpha G)'$, and the existence of a unital embedding of $\cD$ in $A_\omega \cap (A \rtimes_\alpha G)'$ is in fact equivalent to the existence of a unital embedding of $\cD$ into $(A^\alpha)_\omega \cap (A \rtimes_\alpha G)'$. The result follows.
        \end{proof}
    \end{lemma}

    Using the Galois correspondence of Izumi \cite{Izumi02} yields the following.

    \begin{theorem}\label{Galois}
        Let $A$ be a unital, simple, separable C*-algebra and let $G \act^\alpha A$ be an action of a finite group by outer automorphisms. If $A \subseteq A \rtimes_\alpha \cD$ is $\cD$-stable, then there exists an isomorphism $\Phi: A \rtimes_\alpha G \simeq (A \rtimes_\alpha G) \otimes \cD$ such that whenever $C$ is a unital C*-algebra satisfying either
          \begin{enumerate}
            \item $A^\alpha \subseteq C \subseteq A$ or
            \item $A \subseteq C \subseteq A \rtimes_\alpha G$,
          \end{enumerate}
          we have $\Phi(C) = C \otimes \cD$.
            \begin{proof}
              Applying \cite[Corollary 6.6]{Izumi02} gives the following two correspondences:
                \begin{enumerate}
                  \item there is a one-to-one correspondence between subgroups of $G$ with intermediate C*-algebras $A^\alpha \subseteq C \subseteq A$ given by
\begin{equation}
                    H \leftrightarrow A^{\alpha_H};
\end{equation}
                  \item there is a one-to-one correspondence between subgroups of $G$ and intermediate C*-algebras $A \subseteq C \subseteq A \rtimes_\alpha G$ given by
\begin{equation}
                    H \leftrightarrow A \rtimes_{\alpha|_H} H.
\end{equation}
                \end{enumerate}
            In particular, there are only finitely many C*-algebras $C$ between either $A^\alpha \subseteq A$ or $A \subseteq A \rtimes_\alpha G$. As all such lie between the $\cD$-stable inclusion $A^ \alpha \subseteq A \rtimes G$, Theorem \ref{countablymanysubalgebras} yields the desired isomorphism. 
            \end{proof}
      \end{theorem}

  \subsection{(Tracial) Rokhlin properties}\label{ssection:dse-tracial-Rokhlin}

  Here we will restrict ourselves to finite groups for simplicity, although many results hold more generally (see \cite{HirshbergWinter07,HirshbergOrovitz13,GardellaHirshberg18}). 

  \begin{defn}
      Let $A$ be a unital, separable C*-algebra. We say that a finite group action $G \act^\alpha A$ has the Rokhlin property if there are pairwise orthogonal projections $(p_g)_{g \in G} \subseteq A_\omega \cap A'$ summing to $1_{A_\omega}$ such that $(\alpha_g)_\omega(p_h) = p_{gh}$ for $g,h \in G$.
    \end{defn}

    \begin{prop}
      Let $A$ be a unital, separable $\cD$-stable C*-algebra. If $G \act^\alpha A$ is an action of a finite group with the Rokhlin property, then $A^\alpha \subseteq A \rtimes_\alpha G$ is $\cD$-stable.
          \begin{proof}
            This follows from \cite[Theorem 3.3]{HirshbergWinter07}, together with Lemma \ref{lem:fixed-point-in-crossed-product-D-stable}.
          \end{proof}
    \end{prop}

    \begin{defn}
        Let $A$ be a unital, separable C*-algebra. We say that a finite group action $G \act^\alpha A$ has the weak tracial Rokhlin property if for all  $\cF \subseteq A$ finite, $\ee > 0$ and $0 \neq a \in A_+$, there are pairwise orthogonal normalized positive contractions $(e_g)_{g \in G} \subseteq A$ such that
              \begin{enumerate}
                \item $1 - \sum_g e_g \precsim a$;\footnote{For two positive elements $x,y$ in a C*-algebra, we write $x \precsim y$ to mean that $x$ is Cuntz-subequivalent to $y$. That is, there are $(r_n)$ in the C*-algebra such that $r_n^*yr_n \to x$. See \cite[Section 2]{HirshbergOrovitz13}.}
                \item $[e_g,x] \approx_\ee 0$ for all $x \in \cF, g \in G$;
                \item $\alpha_g(e_h) \approx_\ee e_{gh}$ for all $g,h \in G$.
              \end{enumerate}
      \end{defn}

      We note that both Rokhlin and weak tracial Rokhlin actions are necessarily outer.

      \begin{prop}\label{prop:dsi-wtr-preserves-Z}
          Let $A$ be a unital, simple, separable, nuclear, $\cZ$-stable C*-algebra. If $G \act^\alpha A$ is an action of a finite group with the weak tracial Rokhlin property. Then $A^\alpha \subseteq A \rtimes_\alpha G$ is $\cZ$-stable.
            \begin{proof}
              Let $k \in \bN$. By \cite[Theorem 5.6]{HirshbergOrovitz13} $A \rtimes_\alpha G$ is tracially $\cZ$-absorbing, meaning there are tracially large (in the sense of \cite{TomsWhiteWinter15}) c.p.c. order zero maps $\phi: M_k \to (A \rtimes_\alpha G)_\omega \cap (A \rtimes_\alpha G)'$, which can be chosen to be c.p.c. order zero maps $\phi: M_k \to A_\omega \cap (A \rtimes_\alpha G)'$ by the proof of \cite[Lemma 5.5]{HirshbergOrovitz13}. These tracially large c.p.c. order zero maps yield sequences of positive contractions $c_1 = (c_{1n}),\dots,c_k = (c_{kn}) \in A_\omega \cap (A \rtimes_\alpha G)'$ such that if $(e_n) = e := 1 - \sum_i c_i^*c_i$, we have
\begin{equation}
              \lim_{n \to \omega} \max_{\tau \in T(A)} \tau(e_n) = 0, \inf_m \lim _{n \to \omega} \min_{\tau \in T(A)} \tau(c_{1n}^m) > 0 
\end{equation}
              and $c_ic_j^* = \delta_{ij}c_1^2$. By \cite[Proposition 4.11]{GardellaHirshberg18} (which is much more general, applicable to all countable amenable groups), $A \subseteq A \rtimes_\alpha G$ has equivariant property (SI) since $A$ has property (SI).\footnote{A unital, separable, simple, nuclear, $\cZ$-stable C*-algebra has property (SI) as in \cite{MatuiSato12}} Consequently there exists $s \in A_\omega \cap (A \rtimes_\alpha G)'$ such that $s^*s = 1 - \sum_i c_i^*c_i$ and $c_1s = s$. Altogether,
                \begin{itemize}
                  \item $c_1 \geq 0$;
                  \item $c_ic_j^* = \delta_{ij}c_1^2$;
                  \item $s^*s + \sum_i c_i^*c_i = 1$;
                  \item $c_1s = s$. 
                \end{itemize}
                As mentioned in the proof of $(iv) \Rightarrow (i)$ of \cite{MatuiSato12}, $\cZ_{n,n+1}$ is the universal C*-algebra generated by $n+1$ elements satisfying the above four relations (see \cite[Proposition 5.1]{RordamWinter10} and \cite[Proposition 2.1]{Sato10}), and consequently we have a unital *-homomorphism $\cZ_{n,n+1} \to A_\omega \cap (A \rtimes_\alpha G)'$. Therefore $\cZ \into A_\omega \cap (A \rtimes_\alpha G)'$, giving that the desired inclusion is $\cZ$-stable by Lemma \ref{lem:fixed-point-in-crossed-product-D-stable}.
            \end{proof}
        \end{prop}

        \begin{cor}\label{weaktracialRokhlinintermediate}
            Let $A$ be a unital, simple, separable, nuclear, $\cZ$-stable C*-algebra and $G \act^\alpha A$ be an action of a finite group with the weak tracial Rokhlin property. There exists an isomorphism $\Phi: A \rtimes_\alpha G \simeq (A \rtimes_\alpha G) \otimes \cZ$ such that whenever $C$ is a unital C*-algebra satisfying either
              \begin{enumerate}
                \item $A^\alpha \subseteq C \subseteq A$ or
                \item $A \subseteq C \subseteq A \rtimes_\alpha G$,
              \end{enumerate}
              we have $\Phi(C) = C \otimes \cZ$. 
                \begin{proof}
                  This results from combining Proposition \ref{prop:dsi-wtr-preserves-Z} together with Theorem \ref{Galois}, making note that this is an outer action. 
                \end{proof}
          \end{cor}

        \subsection{The diagonal inclusion associated to a group action}\label{ssection:dse-diagonal-inclusions} In the von Neumann setting, a certain diagonal inclusion associated to several automor\hyp{}phisms was considered in \cite{Popa89,Kawamuro99,Burstein10}, and they play a role in subfactor theory. Here we consider a  unital C*-algebraic inclusion of the same form.

    \begin{defn}
      Let $A$ be a C*-algebra, $\alpha_1,\dots,\alpha_n \in \Aut(A)$. The diagonal inclusion associated to $\alpha_1,\dots,\alpha_n$ is
\begin{equation}
      B(\alpha_1,\dots,\alpha_n) = \left\{\bigoplus_{i=1}^n \alpha_i(a) \mid a \in A \right\} \subseteq M_n(A).
\end{equation}
    \end{defn}

    If $G \act^\alpha A$ is an action of a finite group, we'll write
    \begin{equation}
        B(\alpha) = \left\{\bigoplus_{g \in G} \alpha_g(a) \mid a \in A\right\} \subseteq M_{|G|}(A).
    \end{equation}
We note that a diagonal $B(\alpha) \subseteq M_{|G|}(A)$ is unique up to unitary conjugation (by permutation unitaries). As $\cD$-stability of an inclusion is preserved under unitary conjugation, there is no ambiguity in speaking of $\cD$-stability of the inclusion $B(\alpha) \subseteq M_{|G|}(A)$.
    
    \begin{prop}\label{diagonalinclusion}
        Let $G \act^\alpha A$ be an action of a countable discrete group on a unital, separable C*-algebra. If $G = \langle g_1,\dots,g_n\rangle$, then $A \subseteq A \rtimes_\alpha G$ is $\cD$-stable if and only if
\begin{equation}
        B(\id_A,\alpha_{g_1},\dots,\alpha_{g_n}) \subseteq M_{n+1}(A)
\end{equation}
        is $\cD$-stable.
          \begin{proof}
            First suppose that $A \subseteq A \rtimes_\alpha G$ is $\cD$-stable. Let $\cF \subseteq \cD, \cG \subseteq M_{n+1}(A)$ be finite and $\ee > 0$. Let $\cG' \subseteq A$ be the set of matrix coefficients of elements of $\cG$, together with the identity of $A$, and let $L := \max\{1,\max_{a \in \cG'}\|a\|\}$. Relabel $\id_A,\alpha_{g_1},\dots,\alpha_{g_n}$ as $\alpha_1,\dots,\alpha_{n+1}$. Let
\begin{equation}
            \delta := \frac{\ee}{(4L+1)(n+1)^2}
\end{equation}
            and let $\psi: \cD \to A$ be a u.c.p. $(\cF,\delta,\cG' \cup\{u_{g_i}\}_{i=1}^n)$-embedding, where $(u_g)$ are the implementing unitaries for $\alpha$. Let $\phi: D \to B(\alpha) \subseteq M_{|G|}(A)$ be given by
            \begin{equation}
            \phi(d) := \bigoplus_{i=1}^{n+1} (\alpha_i \circ \psi)(d).
            \end{equation}
            Clearly $\phi$ will be $(\cF,\delta)$-multiplicative since each component is the composition of a *-homomorphism (which are contractive) with a map which is $(\cF,\delta)$-multiplicative. Now for $d \in \cF$ and $a = (a_{ij}) \in \cG$, we have
\begin{equation}
\begin{split}                \|[\phi(d),(a_{ij})]\| &\leq \sum_{i,j=1}^{n+1}\|\alpha_i(\psi(d))a_{ij} - a_{ij}\alpha_j(\psi(d))\| \\
                                       &\leq \sum_{i,j=1}^{n+1}\| \alpha_i(\psi(d))a_{ij} - \psi(d)a_{ij}\| \\
                                       &\ \ \ + \|\psi(d)a_{ij} - a_{ij}\psi(d)\| + \|a_{ij}\psi(d) - a_{ij}\alpha_j(\psi(d))\| \\
                                       &\leq \sum_{i,j=1}^{n+1}\|a_{ij}\|\left(\|\alpha_i(\psi(d)) - \psi(d)\| + \|\psi(d) - \alpha_j(\psi(d))\|\right) \\
                                       &\ \ \  + \|[\psi(d),a_{i,j}]\| \\
                                       &< (n+1)^2(2L(\delta + \delta) + \delta) \\
                                       &= (n+1)^2(4L+1)\delta = \ee.
\end{split}
\end{equation}

              Conversely if the associated diagonal inclusion is $\cD$-stable we note that if $(x_k) \subseteq B(\id_A,\alpha_{g_1},\dots,\alpha_{g_n})$ is central for $M_{n+1}(A)$, writing
              \begin{equation}
              x_k = \bigoplus_{i=1}^{n+1}\alpha_i(a_k)
              \end{equation}
               yields that $(a_k) \subseteq A$ is central for $A$ and is asymptotically fixed by $\alpha_{g_i}, i=1,\dots,n$. In particular if $\cD \into B(\id_A,\alpha_{g_1},\dots,\alpha_{g_n})_\omega \cap (M_{n+1}(A))'$, then $\cD \into A_\omega \cap (A \rtimes_\alpha G)'$. 
          \end{proof}
      \end{prop}

        \begin{cor}
          Let $G \act^\alpha A$ be an action of a finite group on a unital, separable C*-algebra. Then $A \subseteq A \rtimes_\alpha G$ is $\cD$-stable if and only if
\begin{equation}
          B(\alpha) \subseteq M_{|G|}(A)
\end{equation}
          is $\cD$-stable. 
        \end{cor}

  \section{Examples}\label{section:dse-examples}

  \subsection{Non-examples}\label{ssection:dse-non-examples}
  We first start with some non-examples. Villadsen's C*-algebras with perforation will be useful (see \cite{TomsWinter09} for good exposition). Let $\cQ = \bigotimes_n M_n$ denote the universal UHF C*-algebra.

  \begin{theorem}[\cite{Villadsen98,Toms08}]
      There exists a unital, simple, separable, nuclear C*-algebra $C$ satisfying the UCT such that $C \not\simeq C \otimes \cZ$ and $C$ contains the universal UHF algebra unitally. Moreover $C$ is tracial and can be chosen to be AH with
\begin{equation}
       (K_0(A),K_0(A)^+,[1]_0,K_1(A)) = (\bQ,\bQ_+,1,0).
\end{equation}
    \end{theorem}

      \begin{cor}
        There exists an embedding $\cQ \into \cQ$ which is not $\cZ$-stable. In particular, it is not $\cQ$-stable.
          \begin{proof}
            Let $C$ be as above. Note that $\cQ \subseteq C$ so we must find an embedding $C \into \cQ$. As $C$ is unital, separable, exact, satisfies the UCT and has a faithful amenable trace (it has traces, and every such trace will be faithful and amenable since $C$ is nuclear and simple) and there is clearly a morphism between $K_0$-groups, \cite[Theorem D]{Schafhauser20} gives an embedding $C \into \cQ$. Conse\hyp{}quently there is an embedding
\begin{equation}
            \cQ \into C \into \cQ
\end{equation}
            which is not $\cQ$-stable since there is an intermediate C*-algebra $C$ with $C \not\simeq C \otimes \cZ$. 
          \end{proof}
      \end{cor}

        \begin{cor}
          There is an embedding $\cZ \into \cQ$ which is not $\cZ$-stable. 
            \begin{proof}
              Take $C$ as above and take the chain of embeddings (noting that $\cQ$ is $\cZ$-stable)
\begin{equation}
              \cZ \into \cQ \otimes \cZ \simeq \cQ \into C \into \cQ.
\end{equation}
            \end{proof}
        \end{cor}

  \begin{cor}
    There is an embedding $\cZ \into \cO_2$ which is not $\cZ$-stable. 
      \begin{proof}
        Just take the same embedding as above together with an embedding $\cQ \into \cO_2$. 
      \end{proof}
  \end{cor}

    \begin{remark}
      All *-homomorphisms between strongly self-absorbing C*-alge\hyp{}bras are approximately unitarily equivalent by \cite[Corollary 1.12]{TomsWinter07}, or even asymptotically unitarily equivalent by \cite[Theorem 2.2]{DadarlatWinter09}. Therefore $\cD$-stability is not closed under these equivalences (nor homotopy, see \cite[Corollary 3.1]{DadarlatWinter09}).
    \end{remark} 

    The only method we have used to show that an inclusion is not $\cD$-stable is by finding an intermediate algebra which is not $\cD$-stable.
    There are plenty of examples of stably finite C*-algebras with perforation or higher-stable rank (in particular non-$\cZ$-stable C*-algebras \cite{Rordam04}) \cite{Villadsen98, Villadsen99,ElliottVilladsen00,Toms05,HirshbergRordamWinter07,Toms08b,Toms08,TomsWinter09,deLacerdaMortari09,Tikuisis12}. This gives rise to the following two questions.

    \begin{enumerate}
      \item Is there a unital inclusion $B \subseteq A$ of separable C*-algebras such that whenever $C$ is such that $B \subseteq C \subseteq A$, we have $C \simeq C \otimes \cD$ but $B \subseteq A$ is not $\cD$-stable? Is $\cD$-stability equivalent to every intermediate C*-algebra being $\cD$-stable?
      \item To get non-examples we use stably finite C*-algebras with perforation in between sufficiently regular C*-algebras. Is there a way to do this for purely infinite C*-algebras, or is finiteness the only obstruction? Thus we can ask: if $\cD$ is a purely infinite strongly self-absorbing C*-algebra, is every embedding of $\cD$ into itself $\cD$-stable? More specifically, if $B \subseteq A$ is a unital inclusion of simple, separable, purely infinite C*-algebras, is the inclusion $\cO_{\infty}$-stable?
    \end{enumerate}

    Our third question asks if we can get non-examples arising from dynamical systems. 

      \begin{enumerate}
        \item[3.] Is there a unital, separable $\cD$-stable C*-algebra and a (finite) group action $G \act^\alpha A$ such that $A \rtimes_\alpha G$ is $\cD$-stable, but the inclusion is not? One would need $A \rtimes_\alpha G$ to be $\cD$-stable for non-dynamical reasons. 
      \end{enumerate}

  \subsection{Cyclicly permuting tensor powers}\label{ssetcion:dse-cyclic-tensor-powers}

  Here we give a dynamical example to illustrate the discussion in section \ref{section:dse-crossed-products}. In particular, we can look at a consequence of Corollary \ref{autfixedpoint}. 

    \begin{example}
      Let $p,q \in \bN$ be coprime and consider the $q$th tensor power of the UHF algebra $A = M_{p^{\infty}}^{\otimes q}$. Let us examine the action $\bZ_q \act^\sigma A$ given by cyclically permuting the tensors: 
\begin{equation}
      \sigma(a_1 \otimes \cdots \otimes a_q) = a_2 \otimes \cdots \otimes a_q \otimes a_1.
\end{equation}
      One can prove directly or use \cite{HirshbergOrovitz13} or \cite{AGJP22} in order to conclude that this action has the weak tracial Rokhlin property and consequently that $A^{\sigma} \subseteq A \rtimes_\sigma \bZ_q$ is $\cZ$-stable. 

      Alternatively, one can use techniques similar to \cite{HirshbergWinter07} in order to compute the $K$-theory of the fixed point algebra $A^{\sigma}$ to be 
\begin{equation}
      K_0((M_{p^{\infty}}^{\otimes q})^{\sigma}) \simeq \dlim\left(\bZ^q,   \begin{pmatrix}
          p + \frac{p^q - p}{q} & \frac{p^q -p}{q} & \cdots & \frac{p^q - p}{q} \\
          \frac{p^q - p}{q} & p + \frac{p^q - p}{q} & \cdots & \frac{p^q - p}{q} \\
          \vdots & \vdots & \ddots & \vdots \\
          \frac{p^q - p}{q} & \frac{p^q - p}{q} & \cdots & p + \frac{p^q - p}{q}
      \end{pmatrix} \right), 
\end{equation}
      from which one can show that $K_0(A^{\sigma})$ is $p$-divisible.Then using the fact that $K_0(A^{\sigma})$ is $p$-divisible and $A^{\sigma}$ is AF, it follows that  $A^{\sigma}$ is  $M_{p^{\infty}}$-stable. Using Corollary \ref{autfixedpoint}, we then see that $M_{p^{\infty}} \into (A^{\sigma})_\omega \cap A'$. In particular, we have that $A^{\sigma} \subseteq A \rtimes_\sigma \bZ_q$ is $M_{p^{\infty}}$-stable (since clearly if this embedding is fixed by $\bZ_q$, it will commute with the implementing unitaries as well).  
    \end{example}

      \begin{example}
        Following up on the previous example, if we consider the embedding
\begin{equation}
        B := \left\{   \begin{pmatrix}
            x \\ & \sigma(x) \\ & & \ddots \\ & & & \sigma^{q-1}(x)
        \end{pmatrix} \mid x \in M_{p^{\infty}}^{\otimes q} \right\} \subseteq M_q(M_{p^{\infty}}^{\otimes q}) := A,
\end{equation}
        then $B \subseteq A$ is $M_{p^{\infty}}$-stable by Proposition \ref{diagonalinclusion}.
      \end{example}

  \subsection{The canonical inclusion of the CAR algebra in $\cO_2$}\label{ssection:dse-car-in-O2}
    \begin{example}
      Let $\cO_2 = C^*(s_1,s_2)$ be the Cuntz algebra generated by two isometries \cite{Cuntz77}, and consider the inclusion
\begin{equation}
      M_{2^{\infty}} \simeq \ov{\text{span}}\{s_\mu s_\nu^* \mid |\mu| = |\nu|\} \subseteq \cO_2,
\end{equation}
      where for a word $\mu = \{i_1,\dots,i_p\} \in \{1,2\}^p$, $s_\mu = s_{i_1}\cdots s_{i_p}$. This copy of the CAR algebra is precisely the fixed point subalgebra of the gauge action (see \cite{RaeburnBook}). Consider the endomorphism $\lambda: \cO_2 \to \cO_2$ given by
\begin{equation}
      \lambda(x) := s_1xs_1^* + s_2xs_2^*.
\end{equation}
      We note that a sequence $(x_n)$ is $\omega$-asymptotically central for $\cO_2$ if and only if it is $\omega$-asymptotically fixed by $\lambda$. Indeed, if $(x_n)$ is central, then $\|\lambda(x_n) - x_n\| \to^{n \to \omega} 0$ since $[x_n,s_i] \to 0$ for $i=1,2$. On the other hand if $(x_n)$ is asymptotically fixed by $\lambda$, then the inequalities
 \begin{equation}
         \begin{split}
        \|s_ix_n - x_ns_i\| &= \|s_1x_ns_1^*s_i + s_2x_ns_2^*s_i - x_ns_i\| \leq \|\lambda(x_n) - x_n\|\|s_i\| \\
        \|s_i^*x_n - x_ns_i^*\| &= \|s_i^*x_n - s_i^*s_1s_1^* - s_i^*s_2x_ns_2^*\| \leq \|s_i^*\|\|\lambda(x_n) - x_n\|
        \end{split}
 \end{equation}
        imply that $(x_n)$ is asymptotically central. 

        We note that $\lambda|_{M^{2^{\infty}}}$ is the forward tensor shift if we identify $M_{2^{\infty}} = \bigotimes_{\bN} M_2$ (see for example \cite[Section V.4]{DavidsonBook1}). Now \cite{BSKR93} gives an embedding $\xi: M_2 \into (M_{2^{\infty}})$ such that $\lambda_\omega \circ \xi = \xi$. In particular $M_{2^{\infty}} \into (M_{2^{\infty}})_{\omega} \cap \cO_2'$ so that this inclusion is $M_{2^{\infty}}$-stable. 
    \end{example}

    Thinking of $\cO_2$ as the semigroup crossed product $\cO_2 \simeq M_{2^{\infty}} \rtimes_\lambda \bN$ (see \cite{Rordam95,Rordam21}), any intermediate C*-algebra is automatically CAR stable. Consequently each intermediate subalgebra $M_{2^{\infty}} \rtimes d\bN = C^*(M_{2^{\infty}},s_1^d)$ is $M_{2^{\infty}}$-stable. We can do this all concurrently. 

      \begin{cor}
        There exists an isomorphism $\Phi: \cO_2 \simeq \cO_2 \otimes M_{2^{\infty}}$ such that
\begin{equation}
        \Phi(C^*(M_{2^{\infty}},s_1^d)) \simeq C^*(M_{2^{\infty}},s_1^d) \otimes M_{2^{\infty}}
\end{equation}
        for all $d \in \bN$. The same holds if we replace $M_{2^{\infty}}$ by $\cZ$. 
      \end{cor}

  Now let us play with some diagonal inclusions associated to powers of the Bernoulli shift $\lambda$ on $\cO_2$ above.  This will be similar to what was discussed in Section \ref{ssection:dse-diagonal-inclusions}, except we allow endomorphisms.

    \begin{example}
      
  Consider, for $n \in \bN$, the diagonal inclusion
\begin{equation}
  B_n:= \left\{   \begin{pmatrix}
      x \\ & \lambda(x) \\ & & \ddots \\ & & & \lambda^{n-1}(x)
  \end{pmatrix} \mid x \in \cO_2 \right\} \subseteq M_n(\cO_2) =: A_n.
\end{equation}
  Note that both $A_n$ and $B_n$ are isomorphic to $\cO_2$, and in fact this gives an non-trivial inclusion of $\cO_2$ into itself which is $\cO_2$-stable. This is $\cO_2$-stable since a sequence is asymptotically fixed by $\lambda$ if and only if it asymptotically commutes with the algebra. A similar argument to that of Proposition \ref{diagonalinclusion} will yield that this inclusion is $\cO_2$-stable.
    \end{example}

  One can even restrict the diagonal to elements of the CAR algebra $M_{2^{\infty}} \subseteq \cO_2$ sitting as the fixed point subalgebra of the Gauge action as above.
  \begin{example}
  Consider
\begin{equation}
   B_n^{(2)} := \left\{   \begin{pmatrix}
  x \\ & \lambda(x) \\ & & \ddots \\ & & & \lambda^{n-1}(x) \end{pmatrix} \mid x \in M_{2^{\infty}} \right\} \subseteq M_n(\cO_2) = A_n.
\end{equation}
  This is $M_{2^{\infty}}$-stable for the same reasons as above. This gives another inclusion $M_{2^{\infty}} \simeq B_n^{(2)} \subseteq M_n(\cO_2) \simeq \cO_2$ which is CAR-stable.
  \end{example}

  \bibliographystyle{amsalpha}
  \bibliography{biblio}

\end{document}